\documentclass[11pt,a4paper]{amsart}
\usepackage{amssymb,amsmath}
\usepackage{mathpazo}
\usepackage{color}
\headsep=1cm \numberwithin{equation}{section}
\newtheorem{theorem}{Theorem}[section]
\newtheorem{lemma}[theorem]{Lemma}
\newtheorem{proposition}[theorem]{Proposition}
\newtheorem{corollary}[theorem]{Corollary}

\theoremstyle{definition}
\newtheorem{definition}[theorem]{Definition}
\theoremstyle{remark}
\newtheorem{remark}[theorem]{Remark}
\newtheorem{example}[theorem]{Example}

\newtheorem{conjecture}[theorem]{Conjecture}

\newcommand{\im}{\operatorname{im}}
\newcommand{\grade}{\operatorname{grade}}
\newcommand{\pgrade}{\operatorname{p.grade}}
\newcommand{\Kgrade}{\operatorname{K.grade}}
\newcommand{\Egrade}{\operatorname{E.grade}}
\newcommand{\Hgrade}{\operatorname{H.grade}}
\newcommand{\Agrade}{\operatorname{A.grade}}
\newcommand{\cgrade}{\operatorname{c.grade}}
\newcommand{\Cgrade}{\operatorname{\check{C}.grade}}

\newcommand{\Spec}{\operatorname{Spec}}

\newcommand{\rad}{\operatorname{rad}}

\newcommand{\wAss}{\operatorname{wAss}}

\newcommand{\Ht}{\operatorname{ht}}

\newcommand{\V}{\operatorname{V}}

\newcommand{\Ext}{\operatorname{Ext}}
\newcommand{\Supp}{\operatorname{Supp}}

\newcommand{\Hom}{\operatorname{Hom}}

\newcommand{\Ann}{\operatorname{Ann}}

\newcommand{\depth}{\operatorname{depth}}

\newcommand{\lo}{\longrightarrow}
\newcommand{\fm}{\frak{m}}
\newcommand{\fp}{\frak{p}}
\newcommand{\fq}{\frak{q}}
\newcommand{\fa}{\frak{a}}
\newcommand{\fb}{\frak{b}}

\newcommand{\fn}{\frak{n}}

\begin{document}

\author[Asgharzadeh  and Tousi ]{Mohsen Asgharzadeh and Massoud Tousi}

\title[On the notion of Cohen-Macaulayness ...]
{On the notion of Cohen-Macaulayness for non Noetherian  rings}

\address{M. Asgharzadeh, Department of Mathematics, Shahid Beheshti
University, Tehran, Iran-and-School of Mathematics, Institute for
Research in Fundamental Sciences (IPM), P.O. Box 19395-5746, Tehran,
Iran.} \email{asgharzadeh@ipm.ir}
\address{M. Tousi, Department of Mathematics, Shahid Beheshti
University, Tehran, Iran-and-School of Mathematics, Institute for
Research in Fundamental Sciences (IPM), P.O. Box 19395-5746, Tehran,
Iran.} \email{mtousi@ipm.ir}

\subjclass[2000]{13C14, 13C15, 13A50, 13H10.}

\keywords{ Cohen-Macaulay ring, grade of an ideal, height of an
ideal, non-Noetherian
ring, rings of invariants.
}

\begin{abstract}
There exist many characterizations of Noetherian Cohen-Macaulay
rings in the literature. These characterizations do not remain
equivalent if we drop the Noetherian assumption. The aim of this
paper is to provide some comparisons between some of these
characterizations in non Noetherian case. Toward solving a
conjecture posed by Glaz, we give a generalization of the
Hochster-Eagon result on Cohen-Macaulayness of invariant rings, in
the context of non Noetherian rings.
\end{abstract}

\maketitle

\section{Introduction}
Throughout this paper all rings are commutative, associative, with
identity, and all modules are unital. The theory of Cohen-Macaulay
rings is a keystone in commutative algebra. However, the study of
such rings have mostly been restricted to the class of Noetherian
rings. On the other hand, certain families of non Noetherian rings
and modules have achieved a great deal of significance in
commutative algebra. For example, a surprising result of Hochster
indicates that non vanishing of a certain $\check{C}ech$ cohomology
module of the ring of absolute integral closure of a Noetherian
domain implies the Directed Summand Conjecture, see \cite[Theorem
6.1]{Ho2}. While  Noetherian Cohen-Macaulay modules are studied in
several research papers, not so much is known about them in the non
Noetherian case. To the best of our knowledge,  until 1992, there
was not any idea for extending the concept of Cohen-Macaulayness to
non Noetherian rings. In that time Glaz \cite{G3}, considered the
notion of Cohen-Macaulayness  for not Noetherian rings and
conjectured that invariant subrings of certain types of rings would
be Cohen-Macaulay. Two years later, she \cite[Page 219]{G4} defined
an $R$-module $M$ to be Cohen-Macaulay (in the sense of Glaz) if for
each prime ideal $\fp$ of $R$, $\Ht_{M}(\fp)=\pgrade_{R_{\fp}}(\fp
R_{\fp},M_{\fp})$, where $\pgrade_{R_{\fp}}(\fp R_{\fp},M_{\fp})$ is
the polynomial grade of $\fp R_{\fp}$ on $R_{\fp}$-module $M_{\fp}$.
Unfortunately, coherent regular rings are not Cohen-Macaulay with
this definition. Then, in the same paper, Glaz asked how one can
define a non Noetherian notion of Cohen-Macaulayness such that the
definition coincides with the original one in the Noetherian case,
and that coherent regular rings are Cohen-Macaulay, see \cite[Page
220]{G4}. In the following, we collect Glaz's desired properties of
the notion of Cohen-Macaulayness for non Noetherian rings.

\begin{conjecture}\label{con}
Can one find a definition of the notion of non Noetherian
Cohen-Macaulay rings such that it satisfies the following three
conditions.
\begin{enumerate}
\item[(i)] The definition coincides with the original definition in
the Noetherian case.
\item[(ii)] Coherent regular rings are Cohen-Macaulay.
\item[(iii)] For a coherent regular ring $R$ and a group  $G$ of
automorphisms of $R$, assume that there exists a module retraction
map $\rho:R \lo R^G$ and that $R$ is a finitely generated
$R^G$-module. Then $R^G$ is Cohen-Macaulay.
\end{enumerate}
\end{conjecture}
Then, Hamilton \cite{H1}, \cite{H2}, \cite{H3} has introduced the
concept of weak Bourbaki (height) unmixed rings, as a first step
towards non Noetherian Cohen-Macaulay rings. Hamilton \cite{H2}
added the following two more properties that must be satisfied by
non Noetherian Cohen-Macaulay rings.

\begin{enumerate}
\item[(H1)] $R$ is Cohen-Macaulay if and only if $R[X]$ is Cohen-Macaulay.
\item[(H2)] $R$ is Cohen-Macaulay if and only if $R_{\fp}$ is
Cohen-Macaulay for all prime ideals $\fp$ of $R$.
\end{enumerate}

More recently, Hamilton and Marley \cite{HM} introduced a definition
for non Noetherian Cohen-Macaulayness rings. If a ring $R$ satisfies
their definition,  then we say that $R$ is Cohen-Macaulay in the
sense of Hamilton-Marley. They used the theory of $\check{C}ech$
cohomology modules to show that Cohen-Macaulayness in the sense of
Hamilton-Marley satisfies the assertions (i) and (ii) of Conjecture
\ref{con}. Adopt the assumption of Conjecture \ref{con} (iii) and
assume in addition that $\dim R\leq2$ and $G$ is finite such that
its order is a unit in $R$. Then Hamilton and Marley proved the
assertion (iii) of Conjecture \ref{con}. Also, they proved the if
part of (H1) and (H2) by their definition.

Perhaps it is worth pointing out that there are many
characterizations of Noetherian Cohen-Macaulay rings and  modules.
In the non Noetherian case, these are not necessarily equivalent.
All of these characterizations have been chosen as candidates for
definition of non Noetherian Cohen-Macaulay rings, see Definition
\ref{def1}. The aim of the present paper is to provide some
comparisons  between these definitions in not necessarily Noetherian
case. Also, toward solving Conjecture \ref{con}, we will present a
definition of the notion of Cohen-Macaulayness in not necessarily
Noetherian case.

Let $R$ be a ring and $\fa$ an ideal of $R$. The organization
of this paper is as follows.

In Section 2, we deal with the notion  of grade of ideals on
modules. There are many definitions  for the notion of grade of an
ideal of a non Noetherian ring. To make things easier, after
recalling these definitions, for the convenience of the reader, we
collect some of their  properties. For our propose, it seems to be
better to use the Koszul grade. This notion of grade is based on the
work \cite{Ho1}. We denote the Koszul grade of an ideal $\fa$ on an
$R$-module $M$ by $\Kgrade_R(\fa,M)$.

In Section 3, we explore interrelation between different definitions
of non Noetherian Cohen-Macaulay rings. These definitions include
the Glaz and Hamilton-Marley definitions and the notion of weak
Bourbaki unmixed rings. Assume that $\mathcal{A}$ is a non empty
subclass of the class of all ideals of a ring $R$. We give some
connections between preceding modules and modules that
 are Cohen-Macaulay modules in the sense of $\mathcal{A}$ (note that an
$R$-module $M$ is said to be Cohen-Macaulay in the sense of
$\mathcal{A}$, if the equality $\Ht_{M}(\fa)=\Kgrade_R(\fa,M)$ holds
for all ideals $\fa$ in $\mathcal{A}$). These classes of ideals
include the class of all finitely generated ideals, prime ideals,
maximal ideals and the class of all ideals. Our work in this section
is motivated by observing that the inequality
$\Kgrade_R(\fa,M)\leq\Ht_{M}(\fa)$ holds for all ideals $\fa$ of
$R$.

In Section 4, we construct three methods for introducing examples of
non Noetherian rings which are Cohen-Macaulay in the sense of any
definition of Cohen-Macaulayness that appeared in the present paper.
Our first example provides the Cohen-Macaulayness of the polynomial
ring $R[X_1,X_2,\cdots]$, where $R$ is Noetherian and
Cohen-Macaulay. Our second example implies the Cohen-Macaulayness of
absolute integral closure of Noetherian complete local domains of
prime characteristic. Our third example concludes the
Cohen-Macaulayness of the perfect closure of Noetherian regular
local domains of prime characteristic.

In Section 5, we give another definition of Cohen-Macaulayness. We
call it Cohen-Macaulayness in the sense of generalized
Hamilton-Marley, see Definition \ref{def2}. Concerning Conjecture
\ref{con}, we will present the following theorem.
\begin{theorem}\label{main}
The following  assertions hold.\begin{enumerate}
\item[(i)] A Noetherian ring is Cohen-Macaulay with original
definition in Noetherian case if and only if  it is Cohen-Macaulay
in the sense of generalized Hamilton-Marley.
\item[(ii)] Coherent regular rings are Cohen-Macaulay in the sense
of generalized Hamilton-Marley.
\item[(iii)] Let $R$ be a Cohen-Macaulay ring in the sense of
generalized Hamilton-Marley and $G$ a finite group of automorphisms
of $R$ such that the order of $G$ is a unit in $R$. Assume that $R$
is finitely generated as an $R^G$-module. Then $R^G$ is
Cohen-Macaulay in the sense of generalized Hamilton-Marley.
\item[(iv)] Let $R$ be a Noetherian  Cohen-Macaulay ring. Then
the polynomial ring $R[X_1,X_2,\cdots]$ is Cohen-Macaulay in the
sense of generalized Hamilton-Marley.
\item[(v)] If $R_{\fp}$ is Cohen-Macaulay in the sense of
generalized Hamilton-Marley for all prime ideals $\fp$ of $R$, then
$R$ is Cohen-Macaulay in the sense of generalized Hamilton-Marley.
\end{enumerate}
\end{theorem}
After proving Theorem 1.2, we continue our study of the behavior of
rings of invariants of different types of non Noetherian
Cohen-Macaulay rings. In view of Definition 3.1, our list of the
different definitions of Cohen-Macaulayness, includes
Cohen-Macaulayness in the sense of (finitely generated) ideals, weak
Bourbaki (height) unmixed.

\section{Different types of the notion of grade}

In this  section $\fa$ is an ideal of a commutative  ring $R$ and
$M$ an $R$-module. We first give a general discussion on the notion
of grade. There are many definitions for notion of grade of  $\fa$
on $M$. Grade over not necessarily Noetherian rings was first
defined by Barger \cite{B} and Hochster \cite{Ho1}. After them,
Alfonsi \cite{A} combined the grade notions of them into a more
general notion of grade for non Noetherian rings and modules. In
this section, for the convenience of the reader, we collect some of
their properties. To make things easier, we first recall them.

\begin{definition} Let $\fa$ be an ideal of a ring $R$ and $M$ an
$R$-module. Take $\Sigma$ be the family of all finitely generated
subideals $\fb$ of $\fa$. Here, $\inf$ and $\sup$ are formed in
$\mathbb{Z} \cup \{\pm\infty\}$ with the convention that $\inf
\emptyset=+ \infty$ and $\sup \emptyset=- \infty$.

(i) In order to give the definition of Koszul grade when $\fa$ is
finitely generated  by a generating set
$\underline{x}:=x_{1},\cdots, x_{r}$, we first denote the Koszul
complex related to $\underline{x}$ by
$\mathbb{K}_{\bullet}(\underline{x})$. Koszul grade of $\fa$ on $M$
is defined by
$$\Kgrade_R(\fa,M):=\inf\{i \in\mathbb{N}\cup\{0\} | H^{i}(\Hom_R(
\mathbb{K}_{\bullet}(\underline{x}), M)) \neq0\}.$$ Note that by
\cite[Corollary 1.6.22]{BH} and \cite[Proposition 1.6.10 (d)]{BH},
this does not depend on the choice of generating sets of $\fa$. For
an ideal $\fa$ (not necessarily finitely generated), Koszul grade of
$\fa$ on $M$ can be defined by
$$\Kgrade_R(\fa,M):=\sup\{\Kgrade_R(\fb,M):\fb\in\Sigma\}.$$ By using
\cite[Proposition 9.1.2 (f)]{BH}, this definition coincides with the
original definition for finitely generated ideals.

(ii) A finite sequence $\underline{x}:=x_{1},\cdots,x_{r}$ of
elements of $R$ is called weak regular sequence on $M$ if $x_i$ is a
nonzero-divisor on $M/(x_1,\cdots, x_{i-1})M$ for $i=1,\cdots,r$. If
in addition $M\neq (\underline{x})M$, $\underline{x}$ is called
regular sequence on $M$. The classical grade  of $\fa$ on $M$,
denoted by $\cgrade_R(\fa,M)$, is defined to the supremum of the
lengths of all weak regular sequences on $M$ contained in $\fa$.

(iii) (see \cite[Page 149]{N}) The polynomial grade of $\fa$ on M is
defined by
$$\pgrade_R(\fa,M):=\underset{m\rightarrow\infty}{\lim}
\cgrade_{R[t_1, \cdots,t_m]}(\fa R[t_1, \cdots,t_m],R[t_1,,
\cdots,t_m]\otimes_R M).$$

(iv) In the case that $\fa$ is finitely generated by generating set
$\underline{x}:=x_{1}, \cdots, x_{r}$, the $\check{C}ech$ grade of
$\fa$ on $M$ is defined by $\Cgrade_{R}(\fa,M):=\inf\{i\in
\mathbb{N}\cup\{0\}|H_{\underline{x}}^{i}(M)\neq0\}$, where
$H_{\underline{x}}^{i}(M)$ is denoted the $i$-th cohomology of
$\check{C}ech$ complex of $M$ related to $\underline{x}$.
\cite[Proposition 2.7]{HM} implies that $\inf\{i\in
\mathbb{N}\cup\{0\}|H_{\underline{x}}^{i}(M)\neq0\}=\Kgrade_{R}(\fa,M)$.
So $\Cgrade_{R}(\fa,M)$ does not depend on the choice of the
generating sets of $\fa$. For not necessarily finitely generated
ideal $\fa$ the $\check{C}ech$ grade of $\fa$ on $M$ is defined
$$\Cgrade_{R}(\fa,M):=\sup\{\Cgrade_R(\fb,M):\fb\in\Sigma\}.$$ By the
same argument as (i), this is well-defined.

(v) (see \cite{B}) The $\Ext$ grade of $\fa$ on $M$ is defined by
$$\Egrade_{R}(\fa,M):=\inf\{i\in \mathbb{N}\cup\{0\}|\Ext^{i}_{R}(R/\fa,
M)\neq0\}.$$

(vi) The local cohomology grade of $\fa$ on $M$ is defined by
$$\Hgrade_{R}(\fa,M):=\inf\{i\in
\mathbb{N}\cup\{0\}|H_{\fa}^{i}(M):=\underset{n}
{\varinjlim}\Ext^{i}_{R}(R/\fa^{n}, M)\neq0\}.$$

(vii) Let M be a finitely presented R-module and $N$ an $R$-module.
By defining from \cite{A}, $\grade_R(M,N)\geq n$ if and only if for
every finite complex $$\textbf{P}_\bullet: P_n \lo P_{n-1} \lo
\cdots \lo P_0 \lo M \lo 0$$ of finitely generated projective
R-modules $P_i$, there exists a finite complex
$$\textbf{Q}_\bullet: Q_n \lo Q_{n-1} \lo \cdots \lo Q_0  \lo M \lo
0$$ of finitely generated projective modules $Q_j$, and a chain map
$\textbf{P}_\bullet \lo \textbf{Q}_\bullet$ over $M$ such that the
induced maps: $H^i(\Hom_R (\textbf{Q}_\bullet,N)) \lo H^i(\Hom_R
(\textbf{P}_\bullet,N))$ are zero maps for $0 \leq i < n$.
$\grade_R(M,N)$ is equal to the largest integer $n$ for which the
above condition is satisfied. If no such integer $n$ exists we put
$\grade_R(M,N)=+\infty$.

We now recall the definition of $\grade_R(L,.)$ for a general
$R$-module $L$. By definition, $\grade_R(L,N)\geq n$ if for every
$\ell\in L$, $(0:_R\ell)$ contains a finitely generated ideal
$I_{\ell}$ satisfying $\grade_R(R/I_{\ell},N) \geq n$. \cite[Theorem
7.1.10]{G1} implies that, if $L$ is finitely presented, then two
definitions of $\grade_R(L,N)$ coincide. We shall write
$\Agrade_R(\fa,N)$ instead of $\grade_R(R/ \fa,N)$.
\end{definition}

In the next two propositions, we recall some properties and
relations between different types of the notion of grade that
appeared in Definition 2.1. In what follows we will make use them
several times.

\begin{proposition}\label{pro1}Let $\fa$ be an ideal of a ring $R$ and $M$ an
$R$-module. Then the following hold.
\begin{enumerate}
\item[(i)] Let $\underline{y}:=y_1,\cdots,y_t$
be a regular sequence of elements of $\fa$ on $M$. Then
$$\pgrade_{R}(\fa,M)=t+\pgrade_{R}(\fa,\frac{M}
{\underline{y}M}).$$
\item[(ii)] Let $f :R\longrightarrow S$ be a
flat ring homomorphism. Then $$\Kgrade_{R}(\fa,M)\leq
\Kgrade_{S}(\fa S,M\otimes_R S).$$
\item[(iii)] Let $\fa\subseteq\fb$ be a pair of ideals of $R$. Then $\Kgrade_{R}(\fa,M)\leq
\Kgrade_{R}(\fb,M)$.
\item[(iv)](Change of rings) Let $f :R\longrightarrow S$ be a ring homomorphism
and $N$ an $S$-module. Then $\Kgrade_{R}(\fa,N)=\Kgrade_{S}(\fa
S,N)$.
\item[(v)] Let $f :R\longrightarrow S$ be a faithfully flat ring
homomorphism. Then
$$\Kgrade_{R}(\fa,M)=\Kgrade_{S}(\fa S,M\otimes_R S).$$
\item[(vi)]
$\pgrade_{R}(\fa,M)=\pgrade_{R}(\fp,M)$ for some prime ideal $\fp$
containing $\fa$.
\item[(vii)]If $\fa$ is
finitely generated, then
$$\Agrade_R(\fa,M)=\inf\{\Agrade_{R_{\fp}}(\fp
R_{\fp},M_{\fp})|\fp\in \V(\fa)\cap\Supp_R M\}.$$
\end{enumerate}
\end{proposition}

{\bf Proof.} (i) This is Theorem 15 of chapter 5 in \cite{N}.

(ii) First assume that $\fa$ is finitely generated by generating set
$\underline{x}:=x_1,\cdots,x_n$. The symmetry of Koszul cohomology
and Koszul homology says that
$H_{i}(\mathbb{K}_{\bullet}(\underline{x})\otimes_RM))\cong
H^{n-i}(\Hom_R(\mathbb{K}_{\bullet}(\underline{x}),M))$, see
\cite[Proposition 1.6.10 (d)]{BH}. Thus the claim in this case
follows from \cite[Proposition 9.1.2 (c)]{BH}. The desired result
for not necessarily finitely generated ideals follows from the first
case.

(iii) In the case $\fa\subseteq\fb$ is a pair of finitely generated
ideals of $R$, the claim is in \cite[Proposition 9.1.2 (f)]{BH}. The
claim in general case follows from this.

(iv) First assume that $\fa$ is finitely generated by generating set
$\underline{x}$. The claim follows from the isomorphism $\Hom_R(
\mathbb{K}_{\bullet}(\underline{x}), N)\cong\Hom_S(
\mathbb{K}_{\bullet}(\underline{x})\otimes_RS, N)$. Now, assume that
$\fa$ is a general ideal of $R$ (not necessarily finitely
generated). Then, by the former case, we have $\Kgrade_R(\fa,N)\leq
\Kgrade_S(\fa S,N)$. Now, let $\underline{y}$ be a finite sequence
of elements of $\fa S$. Then there exists a finite sequence
$\underline{x}$ of elements of $\fa$ such that
$\underline{y}S\subseteq \underline{x}S$. Again, by the former case,
$$\Kgrade_S(\underline{y}S,N)\leq
\Kgrade_S(\underline{x}S,N)=\Kgrade_R (\underline{x}R,N)\leq
\Kgrade_R(\fa,N).$$ This completes the proof.

(v) This is in \cite[Lemma 7.1.7 (2)]{G1}.

(vi) This is Theorem 16 of chapter 5 in \cite{N}.

(vii) This is in \cite[Theorem 7.1.11]{G1}. $\Box$

\begin{proposition}\label{pro2}
Let $\fa$ be an ideal of a ring $R$ and $M$ an $R$-module. Then the
following hold.
\begin{enumerate}
\item[(i)]$\cgrade_R(\fa,M)\leq \pgrade_{R}(\fa,M)=\Kgrade_{R}
(\fa,M)=\Cgrade_{R}(\fa,M)=\Agrade_R(\fa,M)$.
\item[(ii)] $\Hgrade_{R}(\fa,M)=\Egrade_{R}(\fa,M)$.
\item[(iii)] If $\fa$ is finitely generated, then
$\Egrade_{R}(\fb,M)=\Kgrade_{R}(\fb,M)$.
\end{enumerate}
\end{proposition}

{\bf Proof.} (i) One can deduce easily, from Proposition \ref{pro1}
(i) that $$\cgrade_R(\fa,M)\leq \pgrade_{R}(\fa,M).$$

Assume that $\Sigma$ runs through all finitely generated subideals
$\fb$ of $\fa$. In light of \cite[Theorem 5.11]{N} we see that
$$\pgrade_{R}(\fa,M)=\sup\{\pgrade_R(\fb,M):\fb\in\Sigma\}.$$ In view
of \cite[Proposition 2.7]{HM}, one has
$$\pgrade_{R}(\fb,M)=\Kgrade_{R}(\fb,M)=
\Cgrade_{R}(\fb,M)$$ for all finitely generated ideals $\fb$ of $R$.
This yields such equalities for all ideals $\fa$ of $R$.

On the other hand, equivalency $(1)\Leftrightarrow (4)$  of
\cite[Theorem 7.1.8]{G1}, says that the equality
$\Agrade_{R}(\fb,M)=\Kgrade_{R}(\fb,M)$ holds for all finitely
generated ideals $\fb$ of $R$. By definition, such equality holds
for any ideals if one can shows that
$$\Agrade_{R}(\fa,M)=\sup\{\Agrade_R(\fb,M):\fb\in\Sigma\}.$$ To see
this, first assume that $\Agrade_{R}(\fa,M)\geq n$. Then one can
find a finitely generated subideal $J$ of $(\fa:_R1)=\fa$ satisfying
$\Agrade_R(J,M) \geq n$. So
$$ n\leq\sup\{\Agrade_R(\fb,M):\fb\in\Sigma\}.$$ Conversely, let $n$
be an integer such that $\sup\{\Agrade_R(\fb,M):\fb\in\Sigma\}\geq
n$. Then there is a finitely generated subideal $\fb_0$ of $\fa$
such that $\Agrade_R(\fb_0,M)\geq n$. So, for any $r$ in $R$ we have
$\fb_0\subseteq (\fa:_Rr)$ and $\Agrade_R(\fb_0,M)\geq n$. Hence
$\Agrade_R(\fa,M)\geq n$.

(ii) This follows from \cite[Proposition 5.3.15]{Str}.

(iii) This  is in \cite[Proposition 6.1.6]{Str}. $\Box$

The assumptions and results of Proposition \ref{pro2} are sharp. To
see an example consider the following.

\begin{example}(i) In Proposition \ref{pro2} (iii) the finitely
generated assumption on $\fa$ is really needed. To see this, let
$R:=\mathbb{F}[x_1,\cdots,x_n,\cdots]/(x_1^1,\cdots,x_n^n,\cdots)$,
where $\mathbb{F}$ is a field. Set $\fa:=(x_1,\cdots,x_n,\cdots)$.
Then by \cite[Page 367]{B}, one has $$\Kgrade_R(\fa,R)=0\neq
\Egrade_R(\fa,R).$$

(ii) Adopt the notation of (i) and Assume that $\Sigma$ runs over
all finitely generated subideals $\fb$ of $\fa$. By Proposition
\ref{pro2} (i), one has
$$\Egrade_R(\fb,R)=\Hgrade_R(\fb,R)=\Kgrade_R(\fb,R)=0.$$ Therefore
$$\Egrade_{R}(\fa,M)\neq\sup\{\Egrade_R(\fb,M):\fb\in\Sigma\},$$
and
$$\Hgrade_{R}(\fa,M)\neq\sup\{\Hgrade_R(\fb,M):\fb\in\Sigma\}.$$

(iii) Let $R:=\mathbb{F}[[X,Y]]$, where $\mathbb{F}$ is a field and
set $M:=\bigoplus_{0\neq r\in(X,Y)} R/rR$. By inspection of
\cite[Page 91]{Str}, we find that $\Egrade_R(\fm,M)=1$ and
$\cgrade_R(\fm,M)=0$. This shows that the inequality of Proposition
\ref{pro2} (i) does not equality in general. However, if $M$ is a
finitely generated module over a Noetherian ring $R$, then
\cite[Theorem 1.2.5]{BH} provides that $\cgrade_R(\fa,M)=
\Egrade_{R}(\fa,M)$ for all ideals $\fa$ of $R$ such that $M\neq \fa
M$.
\end{example}

As an easy application of Proposition \ref{pro2} (ii), we give an
elementary proof of a result of Foxby. He proved the following
result as an immediate application of the  New Intersection Theorem
and it has an important role in \cite{Fo}.

\begin{corollary} (see \cite[Corollary 1.5]{Fo}) Let $(A,\fm)$ be
a Noetherian local ring and $C$ an $A$-module satisfies $C\neq\fm
C$. Then $\Egrade_A(\fm,C)\leq(\dim C\leq)\dim A$.
\end{corollary}

{\bf Proof.} Note that $\Kgrade_A(\fm,C)<\infty$, since $C\neq\fm
C$. By Grothendeick's Vanishing Theorem, $H^i_{\fm}(C)=0$ for all
$i>\dim C$. Now, the claim follows by Proposition \ref{pro2} (ii)
and (i). $\Box$

\section{Relations between different definitions of Cohen-Macaulay rings}

There are many characterizations of Noetherian Cohen-Macaulay
modules in the literature. If we apply these characterizations to
non Noetherian modules, then they are not necessarily equivalent.
The aim of this section is to provide some relations between these
definitions, when we apply them to not necessarily Noetherian rings
and modules.\\

\textbf{3.A. The basic definitions.}  In this subsection we recall
some candidates for the notion of Cohen-macaulayness in the context
of non Noetherian rings and modules. In what follows we need the
notion of weakly associated prime ideals of an $R$-module $M$.
Recall that a prime ideal $\fp$ is weakly associated to $M$ if $\fp$
is minimal over $(0 :_R m)$ for some $m\in M $. We denote the set of
weakly associated primes of $M$ by $\wAss_RM$. Also, in order to
give the Hamilton and Marley definition of Cohen-Macaulayness, we
need to recall the following definitions (a) and (b).

\begin{enumerate}
\item[(a)] (\cite[Definition 2.3]{Sch}) Let $\underline{x} = x_{1}, \cdots,
x_{r}$ be a system of elements of $R$. For $m \geq n$ there exists a
chain map $\varphi^{m} _{n} (\underline{x})
:\mathbb{K}_{\bullet}(\underline{x}^{m})\longrightarrow
\mathbb{K}_{\bullet}(\underline{x}^{n})$, which induces by
multiplication of $(\prod x_{i})^{m-n} $.  $\underline{x}$ is called
weak proregular if for each $n>0$ there exists an $m \geq n$ such
that the maps $H_{i}(\varphi^{m} _{n} (\underline{x})) :H_{i}
(\mathbb{K}_{\bullet}(\underline{x}^{m}))\longrightarrow H_{i}
(\mathbb{K}_{\bullet}(\underline{x}^{n}))$ are zero for all $i \geq
1$.

\item[(b)] (\cite[Definition 3.1]{HM}) A sequence
$\underline{x}:=x_1,\cdots,x_\ell$ is called a parameter sequence on
$R$, if (1) $\underline{x}$ is a weak proregular sequence; (2)
$(\underline{x})R \neq R$, and (3) $H^{\ell}_{\underline{x}}
(R)_{\fp} \neq 0$ for all prime ideals $\fp \in \V(\underline{x}R)$.
Also, $\underline{x}$ is called a strong parameter sequence on $R$
if $x_{1},\cdots, x_{i}$ is a parameter sequence on $R$ for all
$1\leq i \leq \ell$.
\end{enumerate}

Now, we are ready to recall the following definitions of the
different types of Cohen-Macaulay rings.

\begin{definition}\label{def1} Let $R$ be a ring and $M$ an
$R$-module. \begin{enumerate}
\item[(i)] (\cite[Definition 4.1]{HM}) $R$ is called
Cohen-Macaulay in the sense of Hamilton-Marley, if each strong
parameter sequence on $R$ becomes a regular sequence on $R$. We
denote this property by \textmd{HM}.

\item[(ii)] (\cite[Page 219]{G4})  $M$ is called
Cohen-Macaulay in the sense of Glaz, if for each prime ideal $\fp$
of $R$, $\Ht_{M}(\fp)=\Kgrade_{R_{\fp}}(\fp R_{\fp},M_{\fp})$ and
denote this by \textmd{Glaz}.

\item[(iii)] (\cite[Definition 1 and 2]{H2}) Let $\fa$ be a finitely
generated ideal of $R$. Set $\mu(\fa)$, the minimal number of
elements of $R$ that need to generate $\fa$. Assume that for each
ideal $\fa$ with the property $\Ht \fa\geq \mu(\fa)$, we have $\min
(\fa)=\wAss_R( R/\fa)$. A ring with such property is called weak
Bourbaki unmixed. We denote this property by \textmd{WB}.

\item[(iv)] Let $\mathcal{A}$ be a non empty subclass of the class of all
ideals of a ring $R$. We say that $M$ is Cohen-Macaulay in the sense
of $\mathcal{A}$, if $\Ht_{M}(\fa)=\Kgrade_R(\fa,M)$ for all ideals
$\fa$ in $\mathcal{A}$. We denote this property by $\mathcal{A}$.
The classes we are interested in are $\Supp_R(M)$,
$\Supp_R(M)\cap\max (R)$, the class of all ideals and the class of
all finitely generated ideals. We denote them respectively by
$\textmd{Spec}$, $\textmd{Max}$, \textmd{ideals} and \textmd{f.g.
ideals}. \end{enumerate}
\end{definition}

This is clear from the above definition that any zero dimensional
ring is Cohen-Macaulay in the sense of each part of Definition
\ref{def1}. Also, any one dimensional integral domain is
Cohen-Macaulay in
the sense of each part of Definition \ref{def1}.\\

\textbf{3.B. Relations.} The following diagram illustrates our work
in this subsection:\[\begin{array}{lllllll}
\textmd{Max}\Leftarrow\textmd{Spec}
\Leftrightarrow\textmd{ideals}\Rightarrow\textmd{Glaz}
\Rightarrow\textmd{f.g.
ideals}\Rightarrow\textmd{HM}\Leftarrow\textmd{WB} \  \ (\ast).
\end{array}\]
Also, when  the base ring is coherent, we show that
$\textmd{Spec}\Rightarrow\textmd{WB}$.\\

The key to the work in this subsection is given by the following
elementary lemma.

\begin{lemma}\label{key}
Let $\fa$ be an ideal of a ring $R$ and $M$ a finitely generated
$R$-module. Then $\Kgrade_R (\fa,M)\leq \Ht_{M}(\fa)$.
\end{lemma}

{\bf Proof.} If $M/ \fa M=0$, then $\Ht_{M}(\fa)=+\infty$.
Therefore, we can assume that $\Supp_R(\frac{M}{\fa
M})=V(\fa)\cap\Supp M \neq \emptyset$. Let $\fq\in V(\fa)\cap \Supp
M$.  By parts (ii) and (iii) of Proposition 2.2, one gets
$$\Kgrade_R(\fa,M)\leq \Kgrade_{R_{\fq}}(\fa R_{\fq},M_{\fq})\leq
\Kgrade_{R_{\fq}}(\fq R_{\fq},M_{\fq}).$$ Thus, it is enough for us
to show that if $(R,\fm)$ is a quasi local ring and $M$ a finitely
generated non-zero $R$-module, then $\Kgrade_R(\fm,M)\leq \dim M$.
Applying Proposition \ref{pro1} (iv) for the change ring
$R\longrightarrow R/ \Ann M$, we may assume that $M$ is a faithful
$R$-module. So, $\dim M=\dim R$. If $\dim R=\infty$, we have nothing
to prove. Hence we can assume that $\dim R<\infty$.
\cite[Proposition 2.4]{HM} says that $H^i_{\underline{y}}(M)=0$ for
all $i>\dim R =\dim M$ and all finite sequences $\underline{y}$ of
elements of $R$. On the other hand for a finite sequence
$\underline{x}$ of elements of $\fm$, by Nakayama's Lemma,
$M/\underline{x}M\neq 0$, and so
$\Kgrade_R(\underline{x},M)<\infty$. Consequently, by using
Proposition \ref{pro2} (i),
$\Kgrade_R(\fm,M)=\Cgrade_{R}(\fm,M)\leq \dim M$. $\Box$\\

The next result gives the proof of the following implications:
$$\textmd{Spec}
\Leftrightarrow\textmd{ideals}\Rightarrow\textmd{Glaz}
\Rightarrow\textmd{f.g. ideals}.$$

\begin{theorem}Let $M$ be a finitely generated $R$-module.
Consider the following conditions:
\begin{enumerate}
\item[(i)]$\Ht_{M}(\fp)=\Kgrade_R(\fp,M)$ for all prime ideals $\fp$ in $\Supp_ R(M)$.
\item[(ii)]
$\Ht_{M}(\fa)=\Kgrade_R(\fa,M)$ for all ideals $\fa$ of $R$.
\item[(iii)] $\Ht_{M}(\fq)=\Kgrade_{R_{\fp}}(\fq R_{\fp},M_{\fp})$
for all prime ideals $\fp,\fq$ in $\Supp_ R(M)$ with
$\fq\subseteq\fp$.
\item[(iv)] $\Ht_{M}(\fp)=\Kgrade_{R_{\fp}}(\fp R_{\fp},M_{\fp})$
for all prime ideals $\fp$ in $\Supp_ R(M)$.

\item[(v)] $\Ht_{M}(\fa)=\Kgrade_R(\fa,M)$  for all finitely
generated ideals $\fa$ of $R$.
\end{enumerate}
Then $(i)\Leftrightarrow (ii)\Rightarrow (iii)\Rightarrow
(iv)\Rightarrow (v)$.
\end{theorem}

{\bf Proof.} $(i)\Rightarrow (ii)$ Let $\fa$ be an ideal of $R$. By
Proposition \ref{pro1} (vi)  and Proposition \ref{pro2} (i), there
exists a prime ideal $\fp$ of $R$ containing $\fa$ such that
$\Kgrade_R(\fa,M)=\Kgrade_R(\fp,M)$. In view of Lemma \ref{key}, one
can find that$$\Kgrade_R(\fa,M)=\Kgrade_R(\fp,M)=\Ht_{M}(\fp)\geq
\Ht_{M}(\fa)\geq \Kgrade_{R}(\fa,M),$$ which completes the proof.

$(ii)\Rightarrow (i)$ This is trivial.

$(ii)\Rightarrow (iii)$ This follows from the following
$$\Kgrade_R(\fq,M)\leq \Kgrade_{R_{\fp}}(\fq
R_{\fp},M_{\fp})\leq\Ht_{M_{\fp}}(\fq R_{\fp})=\Ht_{M}(\fq),$$where
the last inequality follows from Lemma \ref{key}.

$(iii)\Rightarrow (iv)$ This is trivial.

$(iv)\Rightarrow (v)$ Let $\fa$ be a finitely generated ideal of
$R$. Then, Proposition \ref{pro1} (vii), Proposition \ref{pro2} (i)
and our assumption, imply that
\[\begin{array}{ll}
\Kgrade_R(\fa,M)&=\inf\{\Kgrade_{R_{\fp}}(\fp
R_{\fp},M_{\fp})|\fp\in \V(\fa)\cap\Supp_R
M\}\\&=\inf\{\Ht_{M_{\fp}}(\fp R_{\fp})|\fp\in \V(\fa)\cap\Supp_R
M\}\\&=\Ht_M(\fa),
\\
\end{array}\]
which completes the proof. $\Box$

In view of \cite[Proposition 4.10]{HM}, any weak Bourbaki unmixed
ring is  Cohen-Macaulay in the sense of Hamilton-Marley. Thus, in
order to complete the proof of all of desired implications of the
diagram $(\ast)$, we need to state the following.

\begin{theorem}
Let $R$ be a Cohen-Macaulay ring in the sense of finitely generated
ideals. Then $R$ is Cohen-Macaulay in the sense of Hamilton-Marley.
\end{theorem}

{\bf Proof.} Let $\underline{x}:=x_1,\cdots,x_\ell$ be a strongly
parameter sequence on $R$. By equivalency $(a)\Leftrightarrow (c)$
of \cite[Proposition 4.2]{HM}, its enough to show that
$\Kgrade_R(\underline{x}R,R)=\pgrade_R(\underline{x}R,R)=\ell$. For
a finite sequence $\underline{y}:=y_1,\cdots,y_m$ of elements of
$R$,  \cite[Proposition 3.6 ]{HM} state that $\Ht
(\underline{y}R)\geq m$, if $\underline{y}$ is a parameter sequence
on $R$. Now, let $\fq\in \V(\underline{x}R)$ be such that $\Ht
(\fq)= \Ht (\underline{x}R)$. Also, from definition, one has
$$\Kgrade_{R_\fq}(\underline{x} R_\fq,R_\fq)\leq\mu(\underline{x}
R_\fq)\leq\ell.$$ Then, it turns out that
$$\Kgrade_R(\underline{x}R,R)\leq \Kgrade_{R_\fq}(\underline{x}
R_\fq,R_\fq)\leq\ell\leq
\Ht_R(\fq)=\Ht(\underline{x}R)=\Kgrade_R(\underline{x}R,R).$$
Therefore, $\Kgrade_R(\underline{x}R,R)=\ell$, as claimed. $\Box$

Theorem 3.10 is one of our main results in this subsection. To prove
it, we need a couple of lemmas.

\begin{lemma}
Let $R$ be a Cohen-Macaulay ring in the sense of (finitely
generated) ideals and $x$ a regular element of $R$. Then $R/xR$ is
Cohen-Macaulay in the sense of (finitely generated) ideals. In
particular, a ring $A$ is Cohen Macaulay in the sense of (finitely
generated) ideals, if either $A[[X]]$ or $A[X]$ is as well.
\end{lemma}

{\bf Proof.} Let $\fb:=\fa/xR$ be an ideal (resp. finitely generated
ideal) of $R/xR$. By parts $(i)$ and $(iv)$ of Proposition
\ref{pro1}, one can find that
$$\Kgrade_{R/xR}(\fb,R/xR)=\Kgrade_{R}(\fa,R/xR)=
\Kgrade_R(\fa,R)-1.$$ Then it yields that:\[\begin{array}{ll}
\Kgrade_R(\fa,R)-1&= \Kgrade_{R/xR}(\fb,R/xR)\\&\leq\Ht_{R/xR}
(\fb)\\&\leq\Ht_{R} (\fa)-1\\&=\Kgrade_R(\fa,R)-1,
\\
\end{array}\]
which completes the proof. $\Box$

\begin{remark}
(i) There exists an example of a quasi-local ring $R$ such that it
is Cohen–Macaulay in the sense of Hamilton-Marley  but not $R/xR$
for some regular element $x$ of $R$, see \cite[Example 4.9]{HM}.

(ii) Assume that $(R,\fm)$ is a quasi local ring, which is
equidimensional, semicatenary and weak Bourbaki unmixed. Let $x$ be
a regular element of $R$.  \cite[Theorem D]{H3} shows that $R/xR$ is
weak Bourbaki unmixed.
\end{remark}

Recall that a module is coherent if it is finitely generated and
each of its finitely generated submodule is finitely presented. A
ring is coherent if it is coherent as a module over itself.
Noetherian rings are coherent. There are many examples of non
Noetherian coherent rings. For instance, any non Noetherian
valuation domain is a non Noetherian coherent ring.

\begin{lemma}
Let $R$ be a coherent ring and $\underline{x}:=x_1,\cdots,x_\ell$ a
finite sequence of elements of $R$. Then $H^{i}(\Hom_R(
\mathbb{K}_{\bullet}(\underline{x}),R))$ is finitely generated
$R$-module for all $i$.
\end{lemma}

{\bf Proof.} Let $\textbf{F}^\bullet:0\lo F^{0}\lo\cdots\lo
F^i\stackrel{\varphi^i}\lo F^{i+1}\lo\cdots\lo F^{\ell}\lo 0$ be the
Koszul complex of $R$ related to $\underline{x}$. Let $i$ be an
integer between $0$ and $\ell$. By using the exact sequence
$$F^i\lo F^{i+1}\lo\im\varphi^i\lo 0,$$ we find that $\im\varphi^i$
is finitely presented. Consider the exact sequence
$$0\lo\ker\varphi^i\lo F^i\lo\im\varphi^i\lo 0,$$
in which the maps are the natural one. Keep in mind that $R$ is
coherent. Now, \cite[Theorem 2.5.1]{G1} yields that $\ker\varphi^i$
is finitely presented. From this the claim follows. $\Box$

\begin{remark}
(i) The coherent assumption on $R$ in Lemma 3.7 is really needed. To
see an example, let $A$ be a $\mathbb{C}$-algebra generated by all
degree two monomials of
$\mathbb{C}[X_1,X_2,\cdots]:=\bigcup_{n=1}^{\infty}
\mathbb{C}[X_1,\cdots,n]$ and set $R:= A/(X_1X_2)$. We use small
letters to indicate the images in $R$. Then $(0:_Rx_1^2)=(x_2x_i:
i\in\mathbb{N})$ is not finitely generated. So the first  Koszul
homology related to $x_1^2$ is not finitely generated (cf.
\cite[Example 2]{G2}).

(ii) If Koszul (co)homology modules are finitely generated, then one
can see that the vanishing of first Koszul homology implies the
exactness of Koszul complex. But there exists an example which does
not satisfy this, see \cite[Example 2]{K}.
\end{remark}

\begin{lemma}
Let $R$ be a Cohen-Macaulay ring in the sense of ideals. Then
$\wAss_R(R)=\min(R)$, where $\min(R)$ is the set of all minimal
prime ideals of $R$.
\end{lemma}

{\bf Proof.} It is well known that $\min(R)\subseteq \wAss_R(R)$.
Let $\fp\in\wAss_R(R)$.  Then \cite[Lemma 2.8]{HM} state that
$\pgrade_{R_{\fp}}(\fp R_{\fp},R_{\fp})=0$. By applying Proposition
\ref{pro2} (i), one has $\Kgrade_{R_{\fp}}(\fp R_{\fp},R_{\fp})=0$.
The inequality  $\Kgrade_R(\fp,R)\leq \Kgrade_{R_{\fp}}(\fp
R_{\fp},R_{\fp})$ shows that $\Kgrade_R(\fp,R)=0$. Therefore,
$\Ht_{R}(\fp)=0$, i.e., $\fp\in\min (R)$. $\Box$

Now, we are ready in the position to present our next main result.
\begin{theorem}
Let $R$ be a coherent ring. If $R$ is Cohen-Macaulay in the sense of
ideals, then $R$ is weak Bourbaki unmixed.
\end{theorem}

{\bf Proof.} By Theorem 3.3 and \cite[Theorem 2.4.2]{G1}, $R_{\fp}$
is Cohen-Macaulay in the sense of ideals  and it is coherent for all
prime ideals $\fp$ of $R$. Also, if $R_{\fp}$ is weak Bourbaki
unmixed for any $\fp\in \Spec R$, then by \cite[Theorem 3]{H2}, $R$
is weak Bourbaki unmixed.  Thus, we may and do assume that $R$ is
quasi local. Let $\fa$ be a proper finitely generated ideal of $R$
with the property that $\Ht \fa\geq \mu(\fa)$. Then,
$\Kgrade_R(\fa,R)\leq \mu(\fa)\leq\Ht \fa$. So
$\ell:=\Kgrade_R(\fa,R)= \mu(\fa)=\Ht \fa$, since $R$ is
Cohen-Macaulay in the sense of ideals. Let
$\underline{x}:=x_1,\cdots,x_{\ell}$ be a generating set for $\fa$.
Now, we show that $\underline{x}$ is a strong parameter sequence.
Let $1\leq i< \ell$ and set $\fa_i:=(x_1,\cdots,x_i)R$. As the
reader might have guessed, we consider the following long exact
sequence of $R$-modules and $R$-homomorphisms
\[\begin{array}{ll}
\cdots \lo H^{j}(\Hom_R(
\mathbb{K}_{\bullet}(x_1,\cdots,x_{i}),R))\stackrel{x_{i+1}}\lo
H^{j}(\Hom_R( \mathbb{K}_{\bullet}(x_1,\cdots,x_{i}),R))\lo\\
H^{j+1}(\Hom_R( \mathbb{K}_{\bullet}(x_1,\cdots,x_{i+1}),R))\lo
H^{j+1}(\Hom_R(
\mathbb{K}_{\bullet}(x_1,\cdots,x_{i+1}),R))\lo\cdots.\\
\end{array}\] By Lemma 3.7,
$H^{j}(\Hom_R( \mathbb{K}_{\bullet}(x_1,\cdots,x_{i}), R))$ is
finitely generated for all $j$.  Also, $x_{i+1}$ belongs to the
Jacobson radical of $R$. By using of Nakayama's Lemma, one can find
that
$$\Kgrade_R(\fa_i+x_{i+1}R,R)\leq \Kgrade_R(\fa_i,R)+1.$$ An easy
induction shows that
$$\Kgrade_R(\fa_i+(x_{i+1},\cdots,x_{\ell}),R)\leq
\Kgrade_R(\fa_i,R)+(\ell-i).$$On the other hand,
$\Kgrade_R(\fa_i+(x_{i+1},\cdots,x_{\ell}),R)=\ell$. Hence
$\Kgrade_R(\fa_i,R)\geq i$. This implies that $\Kgrade_R(\fa_i,R)=
i$, since $\fa_i$ can be generated by $i$'s elements. And so by
\cite[Proposition 3.3 (e)]{HM}, $x_1,\cdots,x_{i}$ is a parameter
sequence on $R$. Thus, $\underline{x}$ is a strong parameter
sequence on $R$. In view of Theorem 3.4, $R$ is Cohen-Macaulay in
the sense of Hamilton-Marley. Therefore, $\underline{x}$ forms a
weak regular sequence on $R$. So Lemma 3.5 implies that $R/\fa$ is
Cohen-Macaulay in the sense of ideals. Now, let $\fp\in\wAss_R(R/
\fa)$. Then, Lemma 3.9 shows that $\Ht_{R/\fa}(\fp/ \fa)=0$, i.e.,
$\fp\in\min (\fa)$. $\Box$

\textbf{3.C. Examples.} In this subsection, we provide some
counter-examples to show that non of the following implications are
valid:
\[\begin{array}{lllllll}&& \textmd{WB}
\\
&&\Uparrow
\\
&\textmd{f.g. ideals}\Leftarrow\textmd{Max}\Leftrightarrow &
\textmd{HM }\Rightarrow \textmd{f.g. ideals} \  \ (\ast,\ast).
\\
\end{array}\]

One might ask whether the second statement of Theorem 3.3 is true,
if $\Ht_{R}(\fm)=\Kgrade_R(\fm,R)$ for all maximal ideals $\fm$ of
$R$. This, would not be the case, as the next example shows.

\begin{example}
Let $(R,m)$ be a Noetherian local Cohen-Macaulay ring of dimension
$d>1$. Let $X(d-1):=\{\fp\in\Spec R: \Ht\fp\leq d-1\}$. Set
$M_{d-1}:= \underset{\fp\in X(d-1)}\bigoplus R _{\fp}/\fp R _{\fp}$
and consider $S:=R \ltimes M_{d-1}$, the trivial extension of $R$ by
$M_{d-1}$. Then $S$ is a quasi-local ring with the unique maximal
ideal $\fn:=\fm\ltimes M_{d-1}$. By inspection of \cite[Example
2.10]{HM}, we know that $\Kgrade_{S}(\fn,S)=\Ht(\fn)$. Thus, $S$ is
Cohen-Macaulay in the sense of maximal ideals. Again, in light of
\cite[Example 2.10]{HM}, we see that $\Kgrade_{S}(\fa ,S)=0$ for all
ideals $\fa$ of $S$  with the property that $\rad(\fa)\neq\fn$. Now,
take $a$  be in $\fm$ but not in $\bigcup\{\fp:\fp\in\min(R)\}$. One
has $\rad((a,0)S)\neq \fn$ and $\Ht((a,0)S)\neq0$. This yields that
$S$ is not Cohen-Macaulay in the sense of finitely generated ideals.
Also, by \cite[Example 4.3]{HM}, $S$ is not Cohen-Macaulay in the
sense of Hamilton-Marley.
\end{example}

In view of \cite{Ber}, a ring is called regular if every finitely
generated ideal has finite projective dimension. For example,
valuation domains are coherent and regular. So they are
Cohen-Macaulay in the sense of Hamilton-Marley, see \cite[Theorem
4.8]{HM}. Then, the next result completes our list of
counter-examples to the diagram $(\ast,\ast)$.

\begin{proposition}
Let $(R,\fm)$ be a valuation domain. Then, the following are
equivalent.
\begin{enumerate}
\item[(i)] $R$ is Cohen-Macaulay in the sense of ideals.
\item[(ii)] $R$ is Cohen-Macaulay in the sense of prime ideals.
\item[(iii)] $R$ is Cohen-Macaulay in the sense of  Glaz.
\item[(iv)] $R$ is Cohen-Macaulay in the sense of finitely generated ideals.
\item[(v)]$\dim R\leq1$.
\item[(vi)] $R$ is weak Bourbaki unmixed.
\item[(vii)] $R$ is Cohen-Macaulay in the sense of maximal ideals.
\end{enumerate}
\end{proposition}

{\bf Proof.} Without loss of generality we can assume that $R$ is
not a field. Let $\underline{x}$ be a finite sequence  of nonzero
elements of $\fm$. Since $R$ is a  valuation domain, so there is an
element $r$ such that $rR=(\underline{x})R$. Hence
$\Kgrade_R(\underline{x}R,R)\leq1$.  Thus
$\Kgrade_R(\underline{x}R,R)=1$, because $R$ is a domain. Therefore,
we bring the following statement:
\[\begin{array}{lllllll}
\Kgrade(\fa,R)=1 \textit{ for all non-zero proper  ideals } \fa
\textit{ of }  R. \ \ (\star)
\\
\end{array}\]
The assertions $(i)\Leftrightarrow (ii)$, $(ii)\Rightarrow (iii)$
and $(iii)\Rightarrow(iv)$ are hold by Theorem 3.3.

$(iv)\Rightarrow (v)$ For a contradiction assume that $\dim R>1$.
Since the ideals of $R$ are linearly ordered by means of inclusion,
$R$ has only one prime ideal of height one, say $\fp$.  Let $x\in
\fm\setminus\fp$. Then $\Ht(xR)>1$. So in view of $(\star)$, $R$ is
not Cohen-Macaulay in the sense of finitely generated ideals. This
contradiction shows that $\dim R\leq1$.

$(v)\Rightarrow (ii)$ This is obvious.

$(i)\Rightarrow (vi)$ Any finitely generated ideal of a valuation
domain is principal. So valuation domains are coherent. Therefore,
this implication follows by Theorem 3.10.

$(vi)\Rightarrow (v)$ It is enough to show that any valuation domain
of dimension greater than 1 is not weak Bourbaki unmixed. Assume
that $R$ is of that type. Then there is  the chain $0\subsetneqq \fp
\subsetneqq \fq$ of prime ideals of $R$ such that $\Ht(\fp)=1$. Let
$a\in\fp\setminus\{0\}$ and consider the ideal $\fa:=aR$. Since
ideals of $R$ are linearly ordered by means of inclusion,
$\min(\fa)=\{\fp\}$. Assume that $\min (\fa)=\wAss_R(R/\fa)$. Let
$b\in\fq\setminus\fp$. Then $a,b$ is a weak $R$-sequence of length
2, which is a contradiction with $(\star)$. This shows that $\min
(\fa)\neq\wAss_R( R/\fa)$ and consequently $R$ is not weak Bourbaki
unmixed.

$(ii)\Rightarrow(vii)$ is trivial and the remainder implication
$(vii)\Rightarrow(v)$ follows by $(\star)$. $\Box$

\begin{remark}
Let $(R,\fm)$ be an unique factorization valuation domain which is
not a field. By inspection of $(\star)$ in the proof of Proposition
3.12, one has $\dim R=1$, and so $R$ is Cohen-Macaulay in the sense
of each part of Definition \ref{def1}. Indeed, let $\fp$ be a prime
ideal of $R$ with height one. It is enough to show that $R/\fp$ is a
field. One has $\fp=xR$ for some $x$ in $\fp$, because $R$ is an
unique factorization domain. Let $\fb:=\fa/xR$ be a non zero proper
ideal of $R/ xR$, where $\fa$ is an ideal  of $R$.
Then by $ (\star)$ in the proof of Proposition 3.12, we have
$\Kgrade(\fb,R/xR)=1$ and $\Kgrade(\fa,R)=1$. In light of
Proposition 2.2 (i) one has
$\Kgrade_{R/xR}(\fb,R/xR)=\Kgrade_{R}(\fa,R/xR)=
\Kgrade_R(\fa,R)-1=0$. This contradiction shows that $R/xR$ has no
any non zero proper ideal. Therefore, $R/\fp$ is a field as claimed.
\end{remark}

\section{Examples of Cohen-Macaulay rings}

In this section we will construct some examples of non Noetherian
Cohen-Macaulay rings. Our first example provides the
Cohen-Macaulayness of the ring
$R[X_1,X_2,\cdots]:=\bigcup_{i=1}^{\infty}R[X_1,\cdots,X_i]$, when
$R$ is Noetherian and Cohen-Macaulay. Such result gives us that at
least one of the Hamilton's conditions for an appropriate definition
of non Noetherian Cohen-Macaulay ring.

\begin{theorem}
Let $R$ be a Noetherian Cohen-Macaulay ring. Then the ring
$R[X_1,X_2,\cdots]$ is Cohen-Macaulay in the sense of each part of
Definition \ref{def1}.
\end{theorem}

{\bf Proof.} First, we show that $R':=R[X_1,X_2,\cdots]$ is
Cohen-Macaulay in the sense of prime ideals. Let $\fp$ be a prime
ideal of $R'$. We need to show that the equality $\Kgrade_{R'}
(\fp,R')=\Ht_{R'}(\fp)$ holds. For any positive integer $i$, set
$R_i:=R[X_1,\cdots,X_i]$ and consider the prime ideal
$\widetilde{\fp}_i:=\fp\cap R_i$. Then we have the following chain
of subsets of $R'$:$$\widetilde{\fp}_1\subseteq
\widetilde{\fp}_2\subseteq\cdots\subseteq\widetilde{\fp}_i
\subseteq\widetilde{\fp}_{i+1}\subseteq\cdots.$$

Consider the following, only possibility, cases (a) and (b).

\begin{enumerate}
\item[(a)] For infinitely many $i$'s, the condition
$\widetilde{\fp}_{i}R_{i+1}\subsetneqq \widetilde{\fp}_{i+1}$
satisfies.
\item[(b)] Just only for finitely many $i$'s, the condition
$\widetilde{\fp}_{i}R_{i+1}\subsetneqq \widetilde{\fp}_{i+1}$ holds.
\end{enumerate}

In the case (a), for infinitely many $i$'s  the inequality
$\Ht_{R_{i}}(\widetilde{\fp}_{i})<\Ht_{R_{i+1}}(\widetilde{\fp}_{i+1})$
is true, since
$\Ht_{R_{i}}(\widetilde{\fp}_{i})=\Ht_{R_{i+1}}(\widetilde{\fp}_{i}R_{i+1})$.
Then for such $i$'s,  it turns out that
\[\begin{array}{ll}
\Kgrade_{R'} (\widetilde{\fp}_{i}R',R')&=\Kgrade_{R_{i}}
(\widetilde{\fp}_{i},R_{i})\\&=\Ht_{R_{i}}(\widetilde{\fp}_{i})\\
&<\Ht_{R_{i+1}}(\widetilde{\fp}_{i+1})\\&=\Kgrade_{R_{i+1}}
(\widetilde{\fp}_{i+1},R_{i+1})\\&=\Kgrade_{R'}
(\widetilde{\fp}_{i+1}R',R'),
\\
\end{array}\]
where the first equality follows from Proposition \ref{pro1} (v) and
second from the Cohen-Macaulayness of $R_i$. Hence $\Kgrade_{R'}
(\fp, R')=\infty$ and consequently
$\Kgrade_{R'}(\fp,R')=\Ht_{R'}(\fp)$.

In the case (b), there is an integer $k>0$ such that
$\widetilde{\fp}_{k}R_{k+j}=\widetilde{\fp}_{k+j}$ for all $j>0$. So
$$\fp=\bigcup_{i\geq1}
\widetilde{\fp}_{i}=(\widetilde{\fp}_{1}\cup\cdots\cup\widetilde{\fp}_{k})
\cup(\bigcup_{j\geq1} \widetilde{\fp}_{k}R_{k+j}).$$ In particular,
$\fp$ is finitely generated. Let $\{\alpha_1,\cdots,\alpha_{\ell}\}$
be a generating set for $\fp$. Thus, there is a positive integer as
$m$ such that $\alpha_{j}\in R_m$ for all $1\leq j\leq \ell$. One
can see easily that $((\alpha_1,\cdots,\alpha_{\ell})R_m)R'\cap
R_m=(\alpha_1,\cdots,\alpha_{\ell})R_m$, because $R'/R_m$ is a
faithfully flat ring extension. In particular,
$(\alpha_1,\cdots,\alpha_{\ell})R_m$ is a prime ideal of $R_m$. Now,
by \cite[Lemma 4.1]{H1},  $\Ht
_{R_m}((\alpha_1,\cdots,\alpha_{\ell})R_m)=\Ht_{R'}(\fp)$. Therefore
\[\begin{array}{ll}
\Ht_{R'}(\fp)&=\Ht
_{R_m}((\alpha_1,\cdots,\alpha_{\ell})R_m)\\&=\Kgrade_{R_m}
((\alpha_1,\cdots,\alpha_{\ell})R_m, R_m)\\&=\Kgrade_{R'}
((\alpha_1,\cdots,\alpha_{\ell})R', R')\\&=\Kgrade_{R'} (\fp, R').
\\
\end{array}\]
So  $R'=R[X_1,X_2,\cdots]$ is Cohen-Macaulay in the sense of prime
ideals. Due to Theorem 3.3 we know that $R'$ is Cohen-Macaulay in
the sense of ideals. Also, in view of Theorem 3.4, $R'$ is
Cohen-Macaulay in the sense of Hamilton-Marley. By \cite[Corollary
2.3.4]{G1}, $R'$ is coherence. Thus, Theorem 3.10 implies that $R'$
is weak Bourbaki unmixed. $\Box$

\begin{remark}Let $R$ be a Noetherian Cohen-Macaulay ring and $\fa$
a finitely generated ideal of $R[X_1,X_2,\cdots]$ with the property
that $\Ht \fa\geq \mu(\fa)$. Then by \cite[Theorem 4.2]{H1}, all of
the weak associated primes of $\fa$ have the same height, i.e.,
$R[X_1,X_2,\cdots]$ is weak Bourbaki height unmixed. In particular,
$R[X_1,X_2,\cdots]$ is weak Bourbaki unmixed.
\end{remark}

Let $(R,\fm)$ be a Noetherian local domain and let $R^{+}$ be the
integral closure of $R$ in the algebraic closure of its field of
fractions. Theorem 4.5 provides the Cohen-Macaulayness of $R^{+}$.
To deal with this, we establish the following lemma.

\begin{lemma}
Let $f:R\longrightarrow S$ be a flat and integral ring homomorphism.
If $R$ is Cohen-Macaulay in the sense of ideals, then $S$  is also
Cohen-Macaulay in the sense of ideals.
\end{lemma}

{\bf Proof.}  Let $\fq$ be in $\Spec S$ and set $\fp=\fq\cap R$. In
view of Proposition \ref{pro1} (ii), we have
\[\begin{array}{ll} \Ht\fq&\leq\Ht\fp\\&=
\Kgrade_{R}(\fp,R)\\&\leq \Kgrade_{S}(\fp S,S)\\&\leq
\Kgrade_{S}(\fq,S),\\
\end{array}\]and  so Lemma \ref{key}  completes the proof. $\Box$

Note that by \cite[Theorem 4.5]{AH}, $R^+$ is not coherent, when $R$
is of dimension at least $3$ and of positive characteristic. So in
the next result we can not apply Theorem 3.10 for it.

\begin{theorem}
Let $(R,\fm)$ be a Noetherian complete  local domain. Then the
following holds.
\begin{enumerate}
\item[(i)] If $R$ is of prime
characteristic $p$, then $R^{+}$ is Cohen-Macaulay in the sense of
each part of Definition \ref{def1}.
\item[(ii)]
If $\dim R\geq4$ and $R$ is of mixed characteristic, then $R^{+}$ is
not Cohen-Macaulay in the sense of finitely generated ideals.
\item[(iii)] If $\dim R<3$, then $R^{+}$ is Cohen-Macaulay in the sense of
each part of Definition \ref{def1}.
\item[(iv)] If $\dim R\geq3$ and $R$ containing
a field of characteristic $0$, then $R^+$ is not Cohen-Macaulay in
the sense of finitely generated ideals.
\end{enumerate}
\end{theorem}

{\bf Proof.} By Cohen's Structure Theorem there exists a complete
regular local subring $(A,\fm_{A})$ of $R$ such that $R$ is a
finitely generated $A$-module. Recall that $R^{+}=A^+$. Then,
without loss of generality we can assume that $R$ is regular.

(i) First, we show that $R^{+}$ is Cohen-Macaulay in the sense of of
ideals. In view of \cite[Theorem 5.15]{HH}, $R^+$ is a balanced big
Cohen-Macaulay $R$-algebra, i.e., every system of parameters is
regular on $R^{+}$. Over regular local rings, \cite[6.7,
Flatness]{HH} state that any balanced big Cohen-Macaulay module is
flat. Then, Lemma 4.3 yields that $R^{+}$ is Cohen-Macaulay in the
sense of ideals.

Next, we show that $R^{+}$ is weak Bourbaki unmixed.  Let $\fa$ be a
finitely generated ideal of $R^+$ with the property that $\Ht
\fa\geq \mu(\fa)$. Then, $\Kgrade_{R^+}(\fa,R^+)\leq \mu(\fa)\leq\Ht
\fa$. So $n:=\Kgrade_{R^+}(\fa,R^+)= \mu(\fa)=\Ht \fa$, since $R^+$
is Cohen-Macaulay in the sense of ideals. Let $\{a_1,\cdots,a_n\}$
be a generating set for $\fa$. The ring $R^+$ is a direct union of
module finite ring extensions of $R$. Such ring extensions are
Noetherian, local and complete, since $R$ is local and complete. Let
$A$ be one of them, which contains  $R$ and $a_i$ for all $1\leq i
\leq n$. In view of $A^{+}=R^{+}$, we can assume that $a_i\in R$ for
all $1\leq i \leq n$. Set $\fb:=a_1R+\cdots a_nR$. Then $\fb
R^+=\fa$. Because $R^+$ is an integral extension of $R$, we have
$n=\Ht \fa\leq\Ht \fb\leq n$. So $n:= \mu(\fb)=\Ht \fb$. This
implies that $\{a_1,\cdots,a_n\}$ is a part of a system of parameter
for $R$. Keep in mind that $R^+$ is a balanced big Cohen-Macaulay
$R$-algebra. This say's that $\{a_1,\cdots,a_n\}$ is a regular
sequence on $R^+$. It follows from Lemma 3.5 and Lemma 3.9 that
$\wAss_{R^+}(R^+/\fa)=\min (\fa)$.

(ii) For a contradiction assume that $R^{+}$ is Cohen-Macaulay in
the sense of finitely generated ideals. Then by Theorem 3.4, $R^{+}$
is Cohen-Macaulay in the sense of Hamilton-Marley. Also,
\cite[Proposition 3.6]{AH} state that $R^{+}$ is not a balanced big
Cohen-Macaulay algebra for $R$. Thus, there exists a system of
parameters of $R$ as $\underline{x}:=x_1,\cdots,x_{\ell}$ such that
$\underline{x}$ is not regular sequence on $R^+$. For any $1\leq
i\leq \ell$ set $\underline{x}_i:=x_1,\cdots,x_i$. Then
$\Ht(\underline{x}_i R)=i$, because $R$ is Cohen-Macaulay.
\cite[Theorem 19.4]{Mat} says that regular rings are normal. In
particular, going down theorem holds for the integral extension
$R^+/R$. By applying this, one can find that $\Ht(\underline{x}_i
R^+)=i$. So $\Kgrade_{R^+}(\underline{x}_iR^+,R^+)=i$, because
$R^{+}$ is Cohen-Macaulay in the sense of finitely generated ideals.
By using \cite[Proposition 3.3 (e)]{HM}, one can find that
$\underline{x}_i$ is a parameter sequence on $R^+$.  Therefore,
$\underline{x}$ is a strong parameter sequence on $R^+$. Then
$\underline{x}$ is a regular sequence on $R^+$, since $R^{+}$ is
Cohen-Macaulay in the sense of Hamilton-Marley. This is a
contradiction.

(iii) Let $(R,\fm)$ be a Noetherian  local domain of dimension less
than $3$. One can see easily that $R^+$ is a balanced big
Cohen-Macaulay $R$-algebra. Thus by a same reason as (i), $R^{+}$ is
Cohen-Macaulay in the sense of each part of Definition \ref{def1}.

(iv) By our assumptions, one can see that  $R^+$ is not a balanced
big Cohen-Macaulay $R$-algebra, see e.g. \cite[Page 617]{R}. Then by
a same method as (ii), $R^{+}$ is not Cohen-Macaulay in the sense of
finitely generated ideals. $\Box$

Let $R$ be a domain containing a field of characteristic $p > 0$. We
let $R_{\infty}$ denote the perfect closure of $R$, that is,
$R_{\infty}$ is the ring obtained by adjoining to $R$ the $p^n$-th
roots of all its elements. The next result gives the
Cohen-Macaulayness of $R_{\infty}$.

\begin{theorem}
Let $(R,\fm)$ be a Noetherian regular local ring of prime
characteristic $p$. Then $R_{\infty}$ is Cohen-Macaulay in the sense
of each part of Definition \ref{def1}.
\end{theorem}

{\bf Proof.} For each positive integer $n$, set $R_n:=\{x\in
R_\infty|x^{p^n}\in R\}$. By using of \cite[Corollary 8.2.8]{BH},
one can find that the $R$-algebra $R_n$ is flat. Since
$R_{\infty}:=\underset{n}{\varinjlim}R_n$, so $R_{\infty}$ is flat
$R$-algebra. Therefore by Lemma 4.3, $R_{\infty}$ is Cohen-Macaulay
in the sense of ideals.

Let $\fa$ be a finitely generated ideal of $R_\infty$ with the
property that  $\Ht \fa\geq \mu(\fa)$. Then $m:=\mu(\fa)=\Ht
\fa=\Kgrade_{R_\infty}(\fa,R_\infty)$. Let $\{a_1,\cdots,a_m\}$ be a
generating set for $\fa$. There is an integer $\ell$ such that
$a_i\in R_{\ell}$ for all $1\leq i \leq m$. Set
$\fb:=a_1R_{\ell}+\cdots a_mR_{\ell}$. In order to pass from $R$ to
$R_{\ell}$ assume that $R$ is $d$-dimensional. So $\fm$ can be
generated by $d$ elements, namely $x_1,\cdots,x_d$. The ring
$R_{\ell}$ is local with the maximal ideal
$(x_1^{1/p^{\ell}},\cdots,x_d^{1/p^{\ell}})R_{\ell}$. In particular,
$R_{\ell}$ is regular.  Hence we can replace $R$  by $R_{\ell}$.
Also, $\fb R_\infty=\fa$ and $m:= \mu(\fb)=\Ht \fb$. In view of the
equality $\mu(\fb)=\Kgrade_{R}(\fb,R)$ and  by \cite[Exercise
1.2.21]{BH}, one can generated $\fb$ by an $R$-regular sequence
$\underline{b}:=b_1,\cdots,b_m$.  Keep in mind that $R_\infty$ is a
flat $R$-algebra. Then $\underline{b}$ forms a regular sequence on
$R_\infty$. From this, Lemma 3.5 and Lemma 3.9 we get that
$\min(\fa)=\wAss_{R_\infty}(R_\infty/\fa)$. Therefore, $R_\infty$ is
weak Bourbaki unmixed. $\Box$

The argument of the next result involves the concept of Generalized
Principal Ideal Theorem. By definition, a ring $R$ satisfies
\textmd{GPIT} (for Generalized Principal Ideal Theorem) if
$\Ht(\fp)\leq n$ for each prime ideal $\fp$ of $R$ which is minimal
over an $n$-generated ideal of $R$. Rings, with this property are
denoted by \textmd{GPIT}.  For more details on this, see e.g.
\cite{ADEH}. To see an easy example of non \textmd{GPIT} ring, let
$(V,\fm)$ be an infinite dimensional valuation domain. Then, for any
positive integer $n$ one can find an element $x_n$ such that $\Ht
(x_n V)=n$.

\begin{corollary} Let $(R,\fm)$ be a Noetherian local domain of prime
characteristic $p$. Then the following assertions hold.
\begin{enumerate}
\item[(i)] If $R$ is complete, then
$R^+$ is weak Bourbaki height unmixed.
\item[(ii)]If $R$ is regular, then
$R_\infty$ is weak Bourbaki height unmixed.
\end{enumerate}
\end{corollary}

{\bf Proof.} The proof of (ii) is similar as (i). Thus, we give only
the proof of (i). To do this, first note that by \cite[Theorem
3.3]{H1} over \textmd{GPIT},  weak Bourbaki height unmixed follows
by weak Bourbaki unmixed. Thus, in view of Theorem 4.4 (i), the
claim follows by showing that $R^+$ is \textmd{GPIT}. Due to
\cite[Corollary 2.3]{ADEH} we know that any ring which is integral
over a Noetherian domain is \textmd{GPIT}. Therefore $R^+$ is
\textmd{GPIT}. $\Box$

\section{Cohen-Macaulayness of rings of invariants}

Let $R$ be a commutative ring and $G$ a finite group of
automorphisms of $R$. The subring of invariants defined by $R^G:= \{
x\in R : \sigma(x)=x \textit{ for all } \sigma\in G\}$. Assume that
the order of $G$ is a unit in $R$. Then by a famous result of
Hochster and Eagon \cite[Proposition 13]{HE}, we know that if $R$ is
Noetherian and Cohen-Macaulay, then $R^G$ is as well. Our main aim
of the present section can be regarded as a non Noetherian version
of this result. First, we give the proof of Theorem 1.2. To do this,
we need a new definition for the notion of Cohen-Macaulayness for
arbitrary commutative rings as desired in Theorem 1.2.
\begin{definition} \label{def2}
Let $\underline{x}:=x_{1},\cdots,x_{\ell}$ be a finite sequence of
elements  of a ring $R$.

\begin{enumerate}
\item[(i)] For an $R$-module $L$ set
$\mathbb{K}_{\bullet}(\underline{x};L):=\mathbb{K}_{\bullet}(\underline{x})\otimes_RL$.
Recall that for a pair of integers $m \geq n$, there exists a chain
map $$\varphi^{m} _{n} (\underline{x};L)
:\mathbb{K}_{\bullet}(\underline{x}^{m};L)\longrightarrow
\mathbb{K}_{\bullet}(\underline{x}^{n};L)$$ which induces by
multiplication of $(\prod x_{i})^{m-n} $. We call $\underline{x}$ a
generalized proregular sequence on $R$ if for each positive integer
$n$ and any finitely generated $R$-module $M$, there exists an
integer $m \geq n$ such that the maps
$$H_{i}(\varphi^{m} _{n} (\underline{x};M)) :H_{i}
(\mathbb{K}_{\bullet}(\underline{x}^{m};M))\longrightarrow H_{i}
(\mathbb{K}_{\bullet}(\underline{x}^{n};M))$$ are zero for all $i
\geq 1$.
\item[(ii)] We say that $\underline{x}$ is a generalized parameter sequence
on $R$, if (1) $\underline{x}$ is a generalized proregular sequence,
(2) $(\underline{x})R \neq R$, and (3) $H^{\ell}_{\underline{x}}
(R)_{\fp} \neq 0$ for all prime ideals $\fp \in \V(\underline{x}R)$.
\item[(iii)] We call $\underline{x}$ a generalized strong parameter
sequence on $R$, if $x_{1},\cdots, x_{i}$ is a parameter sequence on
$R$ for all $1\leq i \leq \ell$.
\item[(iv)] We say that $R$ is Cohen-Macaulay in the sense of
generalized Hamilton-Marley, if each generalized strong parameter
sequence on $R$ is a regular sequence on $R$.
\end{enumerate}\end{definition}

\begin{remark} (i) Assume that $R$ is a Noetherian ring. Let $\underline{x}
:=x_{1},\cdots,x_{\ell}$ be a finite sequence of elements  of $R$
and $m\geq n$ a pair of positive integers. \cite[ Lemma 4.3.3]{Str}
says that the morphisms $H_{i}(\varphi^{m} _{n} (\underline{x};R))
:H_{i} (\mathbb{K}_{\bullet}(\underline{x}^{m};R))\longrightarrow
H_{i} (\mathbb{K}_{\bullet}(\underline{x}^{n};R))$ are finally null.
Now, let $M$ be a finitely generated $R$-module. By making
straightforward modification of \cite[Lemma 4.3.3]{Str}, one can see
that the following homomorphisms $$H_{i}(\varphi^{m} _{n}
(\underline{x};M)) :H_{i}
(\mathbb{K}_{\bullet}(\underline{x}^{m};M))\longrightarrow H_{i}
(\mathbb{K}_{\bullet}(\underline{x}^{n};M))$$ are finally null. Then
any finite sequence of elements of $R$ is a generalized proregular
sequence.

(ii) Generalized parameter sequence does not coincide with (partial)
systems of parameters if the ring is Noetherian and local. To see an
example, let $\mathbb{F}$ be a field and consider the ring
$R:=\mathbb{F}[[X,Y,Z]]/ (X)\cap(Y,Z)$. We use small letters to
indicate the images in $R$. As was shown by \cite[Theorem 14.1
(ii)]{Mat}, $y$ is a partial systems of parameter. Note that
$\min(yR)=\fp:=(y,z)$, and so $\Ht\fp=0$. By using Grothendieck
Vanishing Theorem, $H^1_y(R)_{\fp}=0$. Therefore, $y$ is not a
generalized parameter sequence.

(iii)  If $(R,\fm)$ is a $d$-dimensional Noetherian local ring, then
by \cite[Theorem 14.1 (ii)]{Mat}, there exists a choice
$\underline{x}:=x_1,\cdots,x_{d}$ of system of parameters such that
$\Ht(x_1,\cdots,x_i)=i$ for all $1\leq i\leq d$. Then $\Ht(\fp)=i$
for all $\fp\in\min(x_1,\cdots,x_i)$ and by applying Grothendieck
non-vanishing theorem, $H^i_{x_1,\cdots,x_i}(R)_{\fp}\neq0$. This
yields that $\underline{x}$ is a generalized strong  parameter
sequence.

(iv) For convention, the ideal generated by the empty sequence is
the zero ideal and the empty sequence is a regular sequence of
length zero over any ring.
\end{remark}

\begin{lemma}
Let $R$ be a  ring. Then the following assertions hold.
\begin{enumerate}
\item[(i)] Assume that $R$ is Noetherian and Cohen-Macaulay.
Then the ring $R[X_1,X_2,\cdots]$ is Cohen-Macaulay in the sense of
generalized Hamilton-Marley.
\item[(ii)] If $R_{\fp}$ is Cohen-Macaulay in the sense of
generalized Hamilton-Marley for all prime ideals $\fp$ of $R$, then
$R$ is Cohen-Macaulay in the sense of generalized Hamilton-Marley.
\end{enumerate}
\end{lemma}

{\bf Proof.} (i) Note that if a ring is Cohen-Macaulay in the sense
of Hamilton-Marley, then it is Cohen-Macaulay in the sense of
generalized Hamilton-Marley.  So (i) follows from Theorem 4.1.

(ii) Let $\underline{x}$ be a generalized strong parameter sequence
on $R$ and $\fp$ a prime ideal containing $\underline{x}$. Let $N$
be a finitely generated $R_{\fp}$-module. One can find a finitely
generated $R$-module as $M$ such that $M_{\fp}\cong N$. Since
$\underline{x}$ is a generalized proregular sequence on $R$ for each
positive integer $n$ there exists an $m \geq n$ such that the maps
$$H_{i}(\varphi^{m} _{n} (\underline{x};M)) :H_{i}
(\mathbb{K}_{\bullet}(\underline{x}^{m};M))\longrightarrow H_{i}
(\mathbb{K}_{\bullet}(\underline{x}^{n};M))$$ are zero for all $i
\geq 1$. On the other hand localization commutes with homology
functors. Therefore, $\underline{x}$ is a generalized proregular
sequence on $R_{\fp}$. By \cite[Proposition 3.3 (c)]{HM},
$\underline{x}$ is a strong parameter sequence on $R_{\fp}$. Hence,
$\underline{x}$ is a generalized strong parameter sequence on
$R_{\fp}$. So, $\underline{x}$ is a regular sequence on $R_{\fp}$
for all prime ideals $\fp$. In particular, $\underline{x}$ is a
regular sequence on $R_{\fp}$ for all prime ideals $\fp$ containing
$\underline{x}R$. Therefore, $\underline{x}$ is a regular sequence
on $R$. $\Box$

The preparation of Theorem 1.2 in the introduction is finished. Now,
we proceed to the proof of it. We repeat Theorem 1.2 to give its
proof.

\begin{theorem}\label{main}
The following  assertions hold.\begin{enumerate}
\item[(i)] A Noetherian ring is Cohen-Macaulay with original
definition in Noetherian case if and only if  it is Cohen-Macaulay
in the sense of generalized Hamilton-Marley.
\item[(ii)] Coherent regular rings are Cohen-Macaulay in the sense
of generalized Hamilton-Marley.
\item[(iii)] Let $R$ be a Cohen-Macaulay ring in the sense of
generalized Hamilton-Marley and $G$ a finite group of automorphisms
of $R$ such that the order of $G$ is a unit in $R$. Assume that $R$
is finitely generated as an $R^G$-module. Then $R^G$ is
Cohen-Macaulay in the sense of generalized Hamilton-Marley.
\item[(iv)] Let $R$ be a Noetherian  Cohen-Macaulay ring. Then
the polynomial ring $R[X_1,X_2,\cdots]$ is Cohen-Macaulay in the
sense of generalized Hamilton-Marley.
\item[(v)] If $R_{\fp}$ is Cohen-Macaulay in the sense of
generalized Hamilton-Marley for all prime ideals $\fp$ of $R$, then
$R$ is Cohen-Macaulay in the sense of generalized Hamilton-Marley.
\end{enumerate}
\end{theorem}

{\bf Proof.} (i) Let $R$ be a Noetherian ring. Note that, in view of
Remark 5.2 (i), any finite sequence of elements of $R$ is a
generalized proregular sequence. First, assume that $R$ is
Cohen-Macaulay with original definition in Noetherian case. Then by
Lemma 5.3 (ii), we may and do assume that $(R,\fm)$ is local. Let
$\underline{x}:=x_1,\cdots,x_{\ell}$ be a strong generalized
parameter sequence for $R$. Due to \cite[Remark 3.2]{HM} we know
that $\Ht(\underline{x}R)=\ell$. In particular, $\underline{x}$ is a
(partial) systems of parameters. So $\underline{x}$ is a regular
sequence on $R$. This shows that $R$ is Cohen-Macaulay in the sense
of generalized Hamilton-Marley.

Now, assume that $R$ is Cohen-Macaulay in the sense of generalized
Hamilton-Marley. Let $\fa$ be an ideal of $R$ of height $\ell$. In
view of \cite[Theorem A.2, Page 412]{BH}, one can find a sequence
$\underline{x}:=x_1,\cdots,x_{\ell}$ of elements of $\fa$ such that
$\Ht(x_1,\cdots,x_i)=i$ for all $1\leq i\leq \ell$ and
$\Ht(\underline{x}R)=\Ht(\fa)$. Then by using \cite[Remark 3.2]{HM},
$\underline{x}$ is a generalized strong parameter sequence on $R$.
Thus, $\underline{x}$ is a regular sequence on $R$, and so
$\cgrade_R(\fa,R)\geq \Ht(\fa)$. Therefore, $R$ is Cohen-Macaulay
with original definition in Noetherian case.

(ii)  \cite[Theorem 4.8]{HM} says that any coherent regular ring is
locally Cohen-Macaulay in the sense of Hamilton-Marley, and so
locally Cohen-Macaulay in the sense of generalized Hamilton-Marley.
Therefore, Lemma 5.3 (ii) implies (ii).

(iii) Let $\underline{x}:=x_1,\cdots,x_\ell$ be a generalized
parameter sequence on $R^G$. In order to show that $\underline{x}$
is a generalized parameter sequence on $R$, we need to show that the
following three assertions hold:

\begin{enumerate}
\item[(a)] $\underline{x}$ is a generalized proregular sequence on $R$,
\item[(b)]$(\underline{x})R \neq R$, and
\item[(c)] $H^{\ell}_{\underline{x}}
(R)_{\fq} \neq 0$ for all prime ideals $\fq \in \V(\underline{x}R)$.
\end{enumerate}

Let $M$  be a finitely generated $R$-module. Since $R$ is a finitely
generated $R^G$-module, we get that $M$ is also finitely generated
as an $R^G$-module. From this one can find easily that
$\underline{x}$ is a generalized proregular sequence on $R$. Hence
(a) is satisfied.

The assertion (b) trivially holds. In order to show (c), assume for
a contradiction that $H^{\ell}_{\underline{x}} (R)_{\fq}= 0$ for
some prime ideal $\fq \in \V(\underline{x}R)$. It follows from
\cite[Page 324, Proposition 23]{Bk} that $S^{-1}(R^G)=(S^{-1}R)^G$
for any multiplicative closed subset $S$ of $R^G$. Set $\fp:=\fq\cap
R^G$ and $S=R^G\setminus\fp$. So $(R^G)_{\fp}\cong(R_{\fp})^G$ and
$\fp\in\V(\underline{x}R^G)$. Since $\underline{x}$ is a parameter
sequence on $R^G$, we have $$0\neq (H^{\ell}_{\underline{x}}
(R^G))_{\fp}\cong H^{\ell}_{\underline{x}}  ((R^G)_{\fp})\cong
H^{\ell}_{\underline{x} } ((R_{\fp})^G).$$ Also,
$H^{\ell}_{\underline{x} } (R_{\fp})_{\fq R_{\fp}}\cong
H^{\ell}_{\underline{x} } (R_{\fq})$. Then, to simplify the
notation, after replacing $R$ by $R_{\fp}$ and $R^G$ by
$(R^G)_{\fp}$, we can assume that $(R^G,\fm)$ is a quasi local ring
with the following properties; $H^{\ell}_{\underline{x}}
(R^G)\neq0$, $\fq\cap R^G=\fm$ and $H^{\ell}_{\underline{x}}
(R)_{\fq}=0$.

Let $\sigma:R \lo R$ be an element of $G$ and $y\in R^G$. Then the
assignment $r/ y^n\mapsto \sigma(r)/ y^n$ induces an $R^G$-algebra
isomorphism $\sigma_y:R_y \lo R_y$. This gives an $R^G$-isomorphism
of the $\check{C}ech$ complexes
$\sigma_1:\Check{\textbf{C}}_{\bullet}(\underline{x},R)\lo
\Check{\textbf{C}}_{\bullet}(\underline{x},R)$. Let $1 \leq i\leq
\ell$. Thus we have an $R^G$-isomorphisms of the $\check{C}ech$
cohomology modules
$$\sigma_2^i:H^i(\Check{\textbf{C}}_{\bullet}(\underline{x},R))\lo
H^i(\Check{\textbf{C}}_{\bullet}(\underline{x},R)).$$ Note that
$\sigma_2^i(tm)=\sigma (t)\sigma_2^i(m)$ for $t\in R$ and $m\in
H^i_{\underline{x}} (R)$. From this one can find that the assignment
$m/s\mapsto \sigma_2^i(m)/\sigma(s)$ for $s\in R\setminus \fq$ and
$m\in H^i_{\underline{x}} (R)$, induces the following
$R^G$-isomorphisms $$\sigma_3^i:H^i_{\underline{x}} (R)_{\fq}\lo
H^i_{\underline{x}} (R)_{\sigma(\fq)}.$$ Assume that $\fq_1$ and
$\fq_2$ are prime ideals of $R$ lying over $\fm$. In view of
\cite[Page 331, Theorem 2 (i)]{Bk}, one can find an element $\sigma$
in $G$ such that $\sigma(\fq_1)=\fq_2$. Also, any maximal ideals of
$R$ contracted to $\fm$. Thus, from the definition of $\sigma_3^i$,
we have $H^{\ell}_{\underline{x}} (R)_{\sigma(\fn)}=0$ for all
$\fn\in\max(R)$ and consequently $H^{\ell}_{\underline{x}} (R)=0$.
Consider the Reynolds operator $\rho:R \lo R^G$. It sends $r\in R$
to $\frac{1}{|G|}\Sigma_{g\in G} gr$. This follows that $R^G$ is a
direct summand of $R$ as $R^G$-module. So $H^{\ell}_{\underline{x}}
(R^G)=0$, a contradiction. This completes the proof of (c).

Now, assume that $\underline{x}$ is a generalized strong parameter
sequence on $R^G$. The same reason as above, shows that
$\underline{x}$ is a  generalized strong parameter sequence on $R$.
Since $R$ is Cohen-Macaulay in the sense of generalized
Hamilton-Marley, we get that $\underline{x}$ is a regular sequence
on $R$. By applying \cite[Proposition 6.4.4 (c)]{BH}, we find that
$\underline{x}$ is a regular sequence on $R^G$. This completes the
proof of (iii).

(iv) and (v) are proved in Lemma 5.3. $\Box$

In the proof of the next result, we use the method of the proof of
Lemma 3.2 (ii) and Lemma 4.1 in \cite{TZ}. Recall that, a group $G$
is said to be locally finite if for every $x\in R$ the orbit of $x$
has finite cardinality.

\begin{lemma}
Let $R$ be a  ring  and $G$ a  group of automorphisms of $R$.
\begin{enumerate}
\item[(i)] Let $\fa$ be an ideal of $R$ and $S$ a pure extension of $R$.
Then $\Kgrade_{R}(\fa, R)\geq \Kgrade_S (\fa S,S)$.
\item[(ii)] Let $\fa$ be an ideal of $R^G$.
Assume that there is a Reynolds operator for the extension $R/R^G$.
Then $\Kgrade_{R^G}(\fa, R^G)\geq \Kgrade_R (\fa R,R)$.
\item[(iii)] Let $\fq$ be a prime ideal of $R$
and $G$ a locally finite group of automorphisms of $R$ such that the
cardinality of orbit of $x$ is a unit in $R$ for every $x\in R$.
Then $\Ht(\fq)\leq\Ht(\fq \cap R^G)$. The equality holds if $G$ is
finite.
\end{enumerate}
\end{lemma}

{\bf Proof.} (i) Let $\underline{y}:=y_1\cdots,y_{s}$ be a finite
sequence of elements of $\fa S$. Then there exists a finite subset
$\underline{x}:=x_1\cdots,x_{\ell}$ of elements of $\fa$ such that
$\underline{y}S\subseteq\underline{x}S$. In view of \cite[Exercise
10.3.31(a)]{BH}, one can find that the natural map $H_i
(\mathbb{K}_{\bullet}(\underline{x}))\lo H_i
(\mathbb{K}_{\bullet}(\underline{x})\otimes_RS)$ is injective for
all $i$. Then, by symmetry of Koszul cohomology and Koszul homology,
one has $\Kgrade_{R}(\underline{x}R, R)\geq
\Kgrade_{R}(\underline{x}R, S)$. Now, by Proposition \ref{pro1}
(iii) and (iv), we find that\[\begin{array}{ll} \Kgrade_{R}(\fa,
R)&\geq \Kgrade_{R}(\underline{x}R, R)\\&\geq
\Kgrade_{R}(\underline{x}R, S)\\&=\Kgrade_{S}(\underline{x}S,
S)\\&\geq \Kgrade_{S}(\underline{y}S, S)
.\\
\end{array}\]So the claim follows from definition.

(ii) By using Reynolds operator, one can find that $R$ is a pure
extension of $R^G$. So (ii) follows from (i).

(iii) Since $G$ is locally finite, so by \cite[Page 323, Proposition
22]{Bk}, the ring extension $R/ R^G$ is integral. The first claim
follows from this. Let
$$\fp_0\subsetneqq\fp_1\subsetneqq\cdots\subsetneqq\fp_n=\fq\cap
R^G$$ be a chain of prime ideals of $R^G$. By lying over theorem,
there exists $\fq_0\in\Spec(R)$ such that $\fq_0\cap R^G=\fp_0$.
Thus by going up theorem, there is a chain of prime ideals of $R$ as
$\fq_0\subsetneqq\fq_1\subsetneqq\cdots\subsetneqq\fq_n$ such that
$\fq_i\cap R^G=\fp_i$. In view of \cite[Page 331, Theorem 2
(i)]{Bk}, there exists an automorphism $\sigma$ in $G$ such that
$\sigma(\fq_n)=\fq$. It is clear that
$$\sigma(\fq_0)\subsetneqq\sigma(\fq_1)\subsetneqq\cdots
\subsetneqq\sigma(\fq_n)=\fq$$ is a chain of prime ideals of $R$ and
so $\Ht\fq\geq\Ht(\fq\cap R^G)$. $\Box$

We now apply Lemma 5.5  to obtain the following  result on the
Cohen-Macaulayness of  rings of invariants in the sense of (finitely
generated) ideals.

\begin{theorem}
Let $R$ be a Cohen-Macaulay ring in the sense of  (finitely
generated) ideals and $G$ a finite group of automorphisms of $R$
such that the order of $G$ is a unit in $R$. Let $\fa$ be a
(finitely generated) ideal of $R^G$. Then $\Kgrade_{R^G}(\fa, R^G)=
\Kgrade_R (\fa R,R)$ and $\Ht(\fa)=\Ht(\fa R)$. In particular, $R^G$
is Cohen-Macaulay in the sense of (finitely generated) ideals.
\end{theorem}

{\bf Proof. } Let $\fa$ be a (finitely generated) ideal of $R^G$ and
$\fq\in \Spec R$ be such that $\Ht(\fa R)=\Ht\fq$. Thus, by Lemma
5.5 (iii), $\Ht(\fa R)=\Ht(\fq\cap R^G)$. Therefore,  Lemma
\ref{key} and Lemma 5.5 (ii) yield that
$$\Ht\fa\geq \Kgrade_{R^G}(\fa, R^G)\geq \Kgrade_R (\fa
R,R)=\Ht(\fa R)=\Ht(\fq\cap R^G)\geq\Ht\fa,$$ which completes the
proof. $\Box$

To complete  our desired list of the behavior of rings of
invariants, on the different types of Cohen-Macaulay rings, we need
to state the following result. A consequence of this is given by
Corollary 5.8.

\begin{proposition}
Let $R$ be a weak Bourbaki (height) unmixed ring and $G$ a finite
group of automorphisms of $R$ such that the order of $G$ is a unit
in $R$. Then $R^G$ is weak Bourbaki (height) unmixed.
\end{proposition}

{\bf Proof.} The proof of weak Bourbaki height unmixed case is
similar as weak Bourbaki unmixed  case. So we give only the proof of
weak Bourbaki unmixed  case. Let $\fa$ be a finitely generated ideal
of $R^G$ with the property that $\Ht \fa\geq \mu(\fa)$. Assume that
$\fp$ belongs to $\wAss_{R^G}(R^G/\fa)$. Then there exists an
element $r$ in $R^G$ such that $\fp\in\min((\fa:_{R^G}r))$. Let
$\fq$ be any prime ideal of $R$ lying over $\fp$. First, we show
that $\fq\in\wAss_{R}(R/\fa R)$. To do this, let $\fq'$ be a prime
ideal of $R$ such that $(\fa R:_{R}r)\subseteq \fq'\subseteq \fq$.
By contraction of this to $R^G$ we get that $\fq'\cap R^G= \fq\cap
R^G$, because $\fa R\cap R^G=\fa$. So $\fq'=\fq$, i.e.,
$\fq\in\wAss_{R}(R/\fa R)$. Let $\fq_0$ be a prime ideal of $R$ such
that $\Ht(\fa R)=\Ht(\fq_0)$. Then, in view of Lemma 5.5 (iii),
$$\Ht(\fa R)=\Ht(\fq_0)=\Ht(\fq_0 \cap R^G)\geq\Ht(\fa) \geq
\mu(\fa) \geq \mu(\fa R).$$This implies that $\fq\in\min(\fa R)$.

Now, we show that $\fp\in\min(\fa)$. To see this, let $\fp'$ be a
prime ideal of $R^G$ and assume that $\fa\subseteq \fp'\subseteq
\fp$.  By lying over theorem, there exists $\fq'\in\Spec(R)$ such
that $\fq'\cap R^G=\fp'$. By applying the going up theorem to this,
we find a prime ideal $\fq''$ of $R$  such that $\fq'\subseteq
\fq''$ and $\fq''\cap R^G=\fp$. As we saw, one has
$\fq''\in\wAss_{R}(R/\fa R)=\min(\fa R)$. This implies that $\fp'=
\fp$ and consequently $\fp\in\min(\fa)$. $\Box$

The statement of the next result involves a non Noetherian version
of the concept of veronese subrings in polynomial ring
$R:=\mathbb{C}[X_1,X_2,\cdots]$. Let $f:=X_{i_1}^{j_1}\cdots
X_{i_{\ell}}^{j_{\ell}}$ be a monomial in $R$. The degree of $f$ is
defined by $d(f):=\sum_{k=1}^{\ell} j_k$. Let $n$ be a positive
integer. We call the $\mathbb{C}$-algebra generated by all monomials
of degree $n$, the $n$-th veronese subring of $R$. We denoted it by
$R_n$.

\begin{corollary}
Let $n$ be a positive integer and let $R_n$ be the $n$-th veronese
subring of $R:=\mathbb{C}[X_1,X_2,\cdots]$. Then $R_n$ is
Cohen-Macaulay in the sense of each part of Definition \ref{def1}.
\end{corollary}

{\bf Proof.} In light of Theorem 4.1 we see that
$\mathbb{C}[X_1,X_2,\cdots]$ is Cohen-Macaulay in the sense of each
part of Definition \ref{def1}. Since $\mathbb{C}$ is an
algebraically closed field, then for each positive integer $n$,
$\mathbb{C}\setminus\{0\}$ has a multiplicative subgroup $G$ of
order $n$.  Let $g$ be in $G$. The assignment $X_i\mapsto g X_i$
induces an action of $G$ on $R$. Assume that $f$ is a monomial in
$R$. Then $f$ belongs to $R^G$ if and only if $g^{d(f)}=1$ for all
$g\in G$. On the other hand by \cite[V. Theorem 5.3]{Ha}, $G$ is
cyclic. So $f$ belongs to $R^G$ if and only if $d(f)=\ell n$ for
some $\ell\in \mathbb{N}\cup\{0\}$. From this we have
$R^G=R_n=\mathbb{C}[f: d(f)\in n\mathbb{N}]$. Due to Theorem 5.6 and
Proposition 5.7 we know that $R_n$ is Cohen-Macaulay in the sense of
ideals and  weak Bourbaki unmixed. Now, the claim follows by Theorem
3.3. $\Box$

It is noteworthy to remark that the converse of the previous results
of this section are not true and their assumptions are really
needed.

\begin{remark}
(i) Let $\mathbb{F}$ be a perfect field of characteristic $2$. In
\cite{Ber}, Bertin presented an action of a finite group $G$ of
order $4$ on $R:=\mathbb{F}[X,Y,Z,W]$ such that $R^G$ is Noetherian
but not Cohen-Macaulay. Thus, in Theorem 5.6 and Proposition 5.7 the
unit assumption on $|G|$ is really needed, even  if $R$ is
Noetherian and regular.

(ii) Let  $A$ be a  Noetherian normal domain which is not
Cohen-Macaulay. In particular, $A$ is a Krull domain. A beautiful
result of Bergman \cite[Proposition 5.2]{Be} state that there is a
principal ideal domain $R$ and an infinite cycle group $G$ such that
$R^G=A$. So, in Theorem 5.6 the finite assumption on $G$ is really
needed, even  if $R^G$ is Noetherian and regular.

(iii) Let $\mathbb{F}$ be a field and set
$R:=\mathbb{F}[[X,Y]]/(XY,Y^2)$. Then $R$ is not Cohen-Macaulay. The
assignments $X\mapsto X$ and $Y\mapsto -Y$ induce an isomorphism
call it $g$. Consider the group of automorphisms generated by $g$
and denote it by $G:=\langle g\rangle$. Then $|G|=2$ and
$R^G=\mathbb{F}[[X]]$, (cf. \cite[Page 448]{F2}). Therefore, the
converse part of Proposition 5.7 is not true, even  if $R^G$ is
Noetherian and regular.

(iv) Fogarty \cite{F2} presented a wild action of a cyclic group $G$
on a local Noetherian ring $R$ such that $ R^G$ is Noetherian and
$\depth R-\depth R^G$ can be arbitrarily large. Thus the assumptions
of $G$ in Lemma 5.5 (ii) is really needed.

(v) Nagata constructed a zero-dimensional Noetherian ring $R$ and a
finite group $G$ of automorphisms of $R$ such that $R^G$ is non
Noetherian, see e.g. the introduction of \cite{F1}. The ring
extension $R/R^G$ is integral, because $G$ is finite. Since $R$ is
zero dimensional, so $R^G$ is zero dimensional. This is clear that
any zero dimensional ring is Cohen-Macaulay in the sense of each
part of Definition \ref{def1}. Thus, $R^G$  is as well. Therefore,
it is possible $R^G$ becomes Cohen-Macaulay without the unit
assumption on $|G|$.
\end{remark}



\begin{thebibliography}{99}

\bibitem[AH]{AH}{I.M. Aberbach}, {M. Hochster}, {\it Finite
tor dimension and failure of coherence in absolute integral
closures}, J.  Pure  Appl. Algebra,  {\bf122}, (1997), 171--184.

\bibitem[A]{A}{B. Alfonsi}, {\it
Grade non Noetherian}, Comm. Algebra, {\bf8}(9), (1981), 811–-840.

\bibitem[ADEH]{ADEH} {D.F. Anderson}, {D.E. Dobbs},
{P.M. Eakin}, {W.J. Heinzer}, {\it On the generalized principal
ideal theorem and Krull domains}, Pacific J. of math., {\bf 146}(2),
(1990), 201--215.


\bibitem[B]{B}{S.F. Barger}, {\it A theory of grade for commutative rings},
Proc. AMS., {\bf36}, (1972), 365–-368.

\bibitem[Be]{Be} {G.M. Bergman}, {\it Groups acting on hereditary rings},
Proc. of London Math. Soc.,{\bf 23}, (1971), 365–-368.


\bibitem[Ber]{Ber} {J. Bertin}, {\it
Anneaux $coh\acute{e}rents\ \ r\acute{e}guliers$}, C. R. Acad. Sci.
Paris, $S\acute{e}r$ A-B, {\bf 273}, (1971).

\bibitem[Bk]{Bk} {N. Bourbaki}, {\it Commutative algebra}, Chapters 1-7,
Springer-Verlag, Berlin, 1989.

\bibitem[BH]{BH} {W. Bruns}, {J. Herzog}, {\it Cohen-Macaulay
rings}, Rev. Ed. Cambridge univ. Press, {\bf 39}, 1998.



\bibitem[F1]{F1} {J. Fogarty}, {\it K\"{a}hler differentials and
Hilbert's fourteenth problem for finite groups}, Amer. J. of math.,
{\bf102}(6), (1980), 1159-1175.


\bibitem[F2]{F2} {J. Fogarty}, {\it On the depth of local rings of
invariants of cyclic groups}, Proc. of AMS., {\bf 83}, (1981),
448--452.

\bibitem[Fo]{Fo} {F.B. Foxby}, {\it On the  $\mu^i$ in a minimal
injective resolution II}, Math. Scan.,  {\bf 41}, (1977), 19-44.


\bibitem[G1]{G1}{S. Glaz}, {\it Commutative coherent rings},
Springer LNM, {\bf1371}, 1989.


\bibitem[G2]{G2}{S. Glaz}, {\it Fixed rings of coherent regular rings},
Comm. Alg., {\bf20}(9), (1992), 2635--2651.

\bibitem[G3]{G3}{S. Glaz}, {\it Coherence, regularity and
homological dimensions of commutative fixed rings}, In: Commutative
algebra,  (Trieste, 1992),  World Sci. Publ., River Edge, NJ,
(1994), 89–-106.

\bibitem[G4]{G4}{S. Glaz}, {\it Homological dimensions of
localizations of polynomial rings}, (Knoxville, TN, 1994), In:
Zero-dimensional commutative rings, In: Lect. Notes in pure and
appl. Math., {\bf 171}, Marcel Dekker, New York, (1995), 209–-222.

\bibitem[Ha]{Ha} {T.W. Haungerford}, {\it Algebra}, Springer
graduate text in math., {\bf73}, 1974.

\bibitem[H1]{H1}{T.D. Hamilton}, {\it Unmixedness and
generalized principal ideal theorem}, Lect. Notes pure appl. Math.,
{\bf241}, (2005), 282--292.

\bibitem[H2]{H2}{T.D. Hamilton}, {\it Weak Bourbaki
unmixed rings: A step towards non Noetherian Cohen-Macaulayness},
Rocky mountain J. of math., {\bf34}(3), (2004), 963–-977.

\bibitem[H3]{H3}{T.D. Hamilton}, {\it Weak Bourbaki unmixed rings:
A step towards non Noetherian Cohen-Macaulayness}, Ph.D. thesis,
University of north Carolina, (1999).

\bibitem[HM]{HM}{T.D. Hamilton}, {T. Marley}, {\it Non Noetherian
Cohen-Macaulay rings}, J. of algebra, {\bf307}, (2007), 343–-360.

\bibitem[Ho1]{Ho1}{M. Hochster}, {\it  Grade-sensitive modules and
perfect modules}, London math. Soc., {\bf29}(3), (1974), 55–-76.

\bibitem[Ho2]{Ho2}{M. Hochster}, {\it Canonical elements in
local cohomology modules and the direct
summand conjecture}, J. of Algebra, {\bf84}, (1983), 503–-553.



\bibitem[HE]{HE}{M. Hochster and J.A. Eagon}, {\it
Cohen-Macaulay rings, invariant theory, and the generic perfection
of determinantal loci}, Amer. J. of math., {\bf 93}, (1971),
1020-1058.

\bibitem[HH]{HH}{M. Hochster}, {C. Huneke}, {\it Infinite integral
extensions and big Cohen-Macaulay algebras}, Ann. of math.,
{\bf135}(2), (1992), 53–-89.

\bibitem[K]{K}{T. Kabele}, {\it  Regularity conditions in non
Noetherian rings}, Trans. AMS., {\bf155}, (1971), 363--374.

\bibitem[Mat]{Mat}{H. Matsumura}, {\it Commutative ring theory},
Cambridge studies in advanced math., {\bf8}, Cambridge, 1989.

\bibitem[N]{N}{D.G. Northcott}, {\it Finite free resolutions},
Cambridge tracts math., {\bf71}, 1976.


\bibitem[R]{R}{ P. Roberts},  {\it Almost Regular Sequences and
the Monomial Conjecture}, Michigan Math. J., {\bf57}, (2008),
615--623.



\bibitem[Sch]{Sch}{P. Schenzel},  {\it
Proregular sequences, local cohomology, and completion}, Math.
Scand., {\bf92}(2), (2003), 271–-289.

\bibitem[Str]{Str}{J.R.  Strooker},  {\it
Homological questions in local algebra}, London math. Lect. Notes
series, {\bf 145}, 1990.

\bibitem[TZ]{TZ}{M. Tousi},  {H. Zakeri}, {\it
Action of certain groups on local cohomology modules and cousin
complexes}, Alg. Colloquium, {\bf 4}(8), (2001), 441--454.
\end{thebibliography}
\end{document}